\documentclass[12pt]{amsart}
\usepackage[utf8]{inputenc}
\usepackage{style} 

\makeatletter
\def\namedlabel#1#2{\begingroup
 #2%
 \def\@currentlabel{#2}%
 \phantomsection\label{#1}\endgroup
}
\makeatother

\title[Hyperbolic Alternating Links in Surfaces with Boundary]{Hyperbolicity of Alternating Links in Thickened Surfaces with Boundary}
\author[C. Adams]{Colin Adams} \address{Department of Mathematics, Williams College, Williamtown, MA 01267}
\email{cadams@williams.edu}
\author[J. Chen]{Joye Chen} \address{Department of Mathematics, MIT,  77 Massachusetts Avenue, Cambridge, MA 02139-4307} \email{joyec@mit.edu}

\begin{document}

\begin{abstract} Let $F$ be a compact orientable surface with nonempty boundary other than a disk. Let $L$ be a link in $F \times I$ with a connected weakly prime cellular alternating projection to $F$. We provide simple conditions that determine exactly when  $(F \times I) \setminus N(L)$ is hyperbolic. We also consider suitable embeddings of $F \times I$ in an ambient manifold $Y$ with boundary and provide conditions on links $L \subset F \times I$ which guarantee hyperbolicity of $Y \setminus N(L)$. These results provide many examples of hyperbolic links in handlebodies and other manifolds. They also provide many examples of  staked links that are hyperbolic.

\end{abstract}

\maketitle

\section{Introduction and Statement of Theorems} 
A link $L$ in a compact 3-manifold $Y$ is called \textit{hyperbolic} if the complement $M = Y \setminus N(L)$ admits a complete metric of constant sectional curvature $-1$. Hyperbolicity has proven useful in studying links in $S^3$, giving rise to many powerful invariants, in particular volume. Hence, the problem of determining which links are hyperbolic is of key interest. Thurston proved that the complement of a link in a compact orientable 3-manifold is hyperbolic if it contains no essential properly embedded spheres, disks, tori or annuli.

In \cite{menasco}, Menasco used this to prove that all non-split prime alternating links in $S^3$ which are not 2-braid links are hyperbolic. This result was extended by Adams et al in 
\cite{small18}, where they proved that all prime, cellular (all complementary regions of the projection are disks) alternating links in thickened closed surfaces of positive genus are hyperbolic. In \cite{hp17}, Howie and Purcell obtained a more general result using angled chunks, proving that under certain conditions, a link $L$ in an arbitrary compact, orientable, irreducible 3-manifold with a weakly prime, cellular alternating projection onto a closed projection surface is hyperbolic. 

A natural next step is to consider projection surfaces with boundary and determine when links which are alternating with respect to these projection surfaces are hyperbolic in manifolds containing these surfaces.  Throughout, we denote a projection surface by $F$, which we require to be connected, orientable, and compact with nonempty boundary. We are interested in links $L \subset F \times I$ and the corresponding 3-manifold $M = (F \times I) \setminus N(L)$, where $N(\cdot)$ denotes a closed regular neighborhood and $I = [0, 1]$ denotes a closed interval. 

Since the manifolds $M$ we are interested in often have higher genus boundary (that is, boundary components with genus at least 2), we would like to have the stronger notion of \textit{tg-hyperbolicity}. 

\begin{definition}
A compact orientable 3-manifold $N$ is \textit{tg-hyperbolic} if, after capping off all spherical boundary components with 3-balls and removing torus boundaries, the resulting manifold admits a complete hyperbolic metric such that all higher genus boundary components are totally geodesic in the metric.
\end{definition}

Ultimately, we would like to use hyperbolic volume to study links in various manifolds, and requiring tg-hyperbolicity allows us to associate a well-defined finite volume to a manifold with higher genus boundary. We also require our links to be prime in $F \times I$ and have cellular alternating projections on $F$.

\begin{definition}
Let $F$ be a projection surface with boundary, and let $L \subset F \times I$ be a link with  projection diagram $\pi(L)$. We say $L$ is \textit{prime} in $F \times I$ if every 2-sphere in $F \times I$ which is punctured twice by $L$ bounds, on one side, a 3-ball intersecting $L$ in precisely one unknotted arc. 

We say $\pi(L)$ is \textit{cellular alternating} on $F$ if it is alternating on $F$ and, after every boundary component of $F$ is capped off with a disk to obtain a closed orientable surface $F_0$ with diagram $\pi(L)$, every complementary region of $F_0 \setminus \pi(L)$ is an open disk.
\end{definition}

When $\pi(L)$ is a \textit{reduced diagram}, there is an easy way to check whether or not $L$ is prime in $F \times I$.

\begin{definition} \label{defn-reduced}
Let $\pi(L)$ be a link projection on a surface $F$ with boundary. We say $\pi(L)$ is \textit{reduced} if there is no circle in $F$ bounding a disk in $F$ and  which intersects $\pi(L)$ transversely in exactly one (double) point. 
\end{definition}

Note that when a projection is not reduced, we can reduce it by flipping that portion of the projection inside the circle and lower the number of crossings. 

\begin{definition} \label{defn-weakly-prime}
Let $\pi(L)$ be a reduced link projection on a surface $F$ with boundary. We say a link projection $\pi(L) \subset F$ is \textit{weakly prime} if every disk $D \subset F$ which has its boundary $\partial D$ intersect $\pi(L)$ transversely in exactly two points contains no crossings of the projection in its interior. 
\end{definition}

We prove the following extensions of Theorem 2 from \cite{small18} and Theorem 1(b) from \cite{menasco} to allow projection surfaces with boundary. 

\begin{proposition} \label{prop-weakly-prime}
Let $F$ be a projection surface with nonempty boundary, and let $L \subset F \times I$ be a link with a connected, reduced, cellular alternating 
projection diagram $\pi(L) \subset F \times \{1/2\}$. Then $L$ is prime in $F \times I$ if and only if $\pi(L)$ is weakly prime on $F \times \{1/2\}$. \newline
\end{proposition}

\subsection{A criterion for hyperbolicity} Our first main result characterizes when alternating links on projection surfaces with boundary are hyperbolic.

\begin{theorem} \label{thm-thickened}
Let $F$ be a projection surface with nonempty boundary which is not a disk, and let $L \subset F \times I$ be a link with a connected, reduced, 
alternating projection diagram $\pi(L) \subset F \times \{1/2\}$ with at least one crossing.  Let $M = (F \times I) \setminus N(L)$. Then $M$ is tg-hyperbolic if and only if the following four conditions are satisfied: 
\begin{enumerate} [label = (\roman*)]
    \item $\pi(L)$ is weakly prime on $F \times \{1/2\}$; 
    \item the interior of every complementary region of $(F \times \{1/2\}) \setminus \pi(L)$ is either an open disk or an open annulus;
    \item if regions $R_1$ and $R_2$ of $(F \times \{1/2\}) \setminus \pi(L)$ share an edge, then at least one is a disk; 
    \item there is no simple closed curve $\alpha$ in $F \times \{1/2\}$ that intersects $\pi(L)$ exactly in a nonempty collection of crossings,  such that for each such crossing, $\alpha$ bisects the crossing and the two opposite complementary regions meeting at that crossing that do not intersect $\alpha$ near that crossing are annuli.
\end{enumerate}

    \end{theorem}

We often refer to $F \times \{1/2\}$ by $F$ when there is no ambiguity. 
We  exclude the case where $F$ is a disk, since this case is covered by \cite{menasco}. Note that in that case, one must also exclude a cycle of bigons. In our case, a cycle of bigons is excluded by conditions (ii) and (iv). Note that condition (ii) implies the link is cellular alternating.

For an alternative formulation, we may start with a closed projection surface $F_0$ with a link $L \subset F_0 \times I$. Choose a set of distinct points $\{x_i\} \in F_0 \setminus \pi(L)$. Then $F := F_0 \setminus \bigcup_{i = 1}^n \mathring{N}(x_i)$ is a surface with boundary and we may consider $L$ as a link in $F \times I$. In this setting, Theorem \ref{thm-thickened} can be rephrased to say that if the $\{x_i\}$ are chosen such that: (i) $\pi(L)$ is weakly prime on $F$, (ii) every region of $F_0 \setminus \pi(L)$ contains at most one of the $x_i$, (iii) no two adjacent regions of $F_0 \setminus \pi(L)$ both contain an $x_i$, and (iv) there is no simple closed curve on $F_0$ intersecting $\pi(L)$ exactly in a nonempty set of crossings such that it bisects each crossing and each of the two regions it does not pass through at each such crossing contain an $x_i$, then $M$ is tg-hyperbolic. Any other choice of $\{x_i\}$ ensures $M$ will not be tg-hyperbolic. 

This formulation fits well with the notion of a staked link, which we discuss in the final section. However, the advantage to the statement of Theorem  \ref{thm-thickened} is it avoids reference to an initial closed surface $F_0$.

\begin{figure}[htbp]
    \centering
    \includegraphics[scale = 0.7]{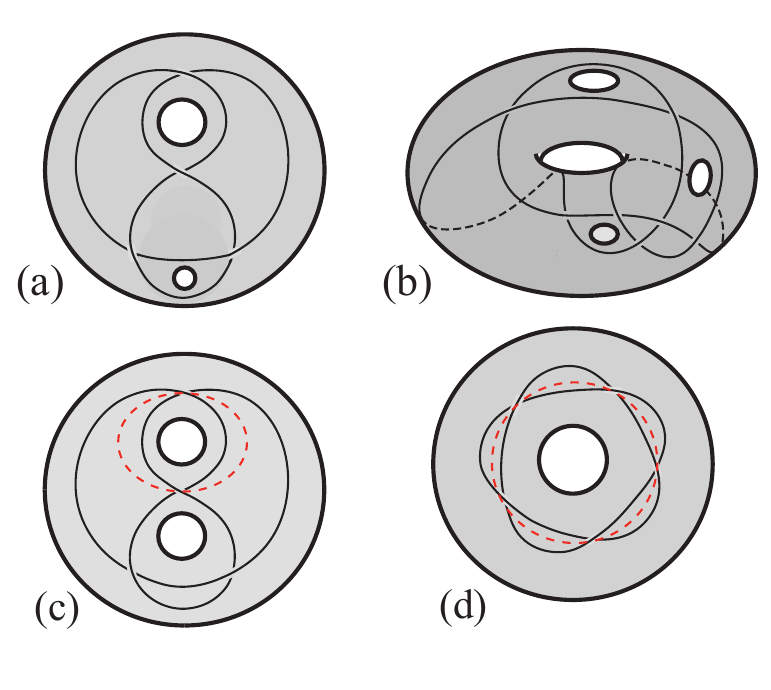}
    \caption{Four examples of $F$ (shaded) and $\pi(L)$ satisfying conditions (i), (ii) and (iii) of Theorem \ref{thm-thickened}.
     Examples (a) and (b) also satisfy condition (iv) so  the corresponding manifolds $(F \times I) \setminus N(L)$  are tg-hyperbolic. Examples (c) and  (d) fail condition (iv) and neither is tg-hyperbolic. A problematic simple closed curve appears in red.} 
    \label{fig-examples-thickened}
\end{figure}

\begin{figure}[htbp]
    \centering
    \includegraphics[scale = 0.7]{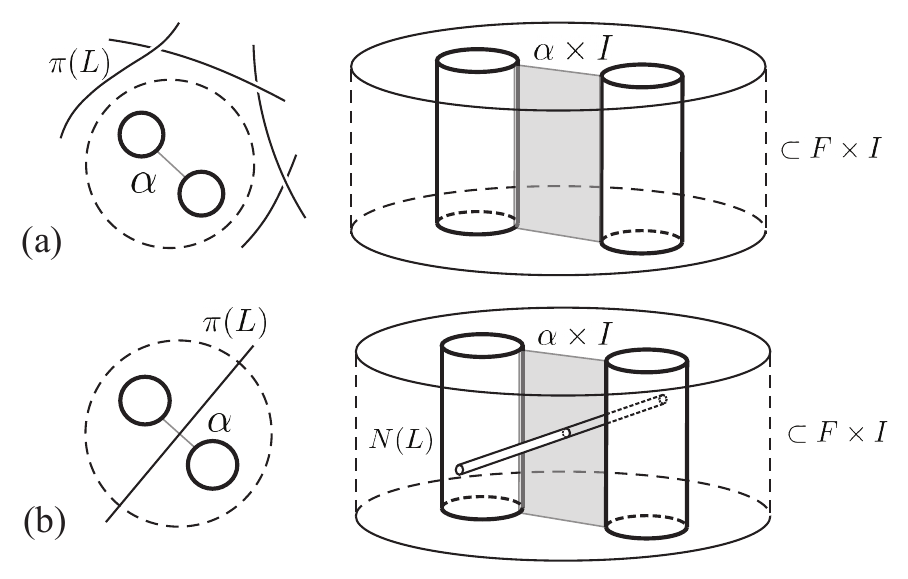}
    \caption{Conditions (ii) and (iii) are necessary for Theorem \ref{thm-thickened}: on the left is a local portion of $\pi(L)$ on $F$ and on the right is the corresponding portion of $M$. In (a) we exhibit an essential disk $\alpha \times I$ when a complementary region has two or more boundary components, and in (b) we exhibit an essential annulus $(\alpha \times I) \setminus N(L)$ when adjacent regions both have boundary.} 
    \label{fig-nonexamples-thickened}
\end{figure}

\begin{figure}[htbp]
    \centering
    \includegraphics[scale = 0.35]{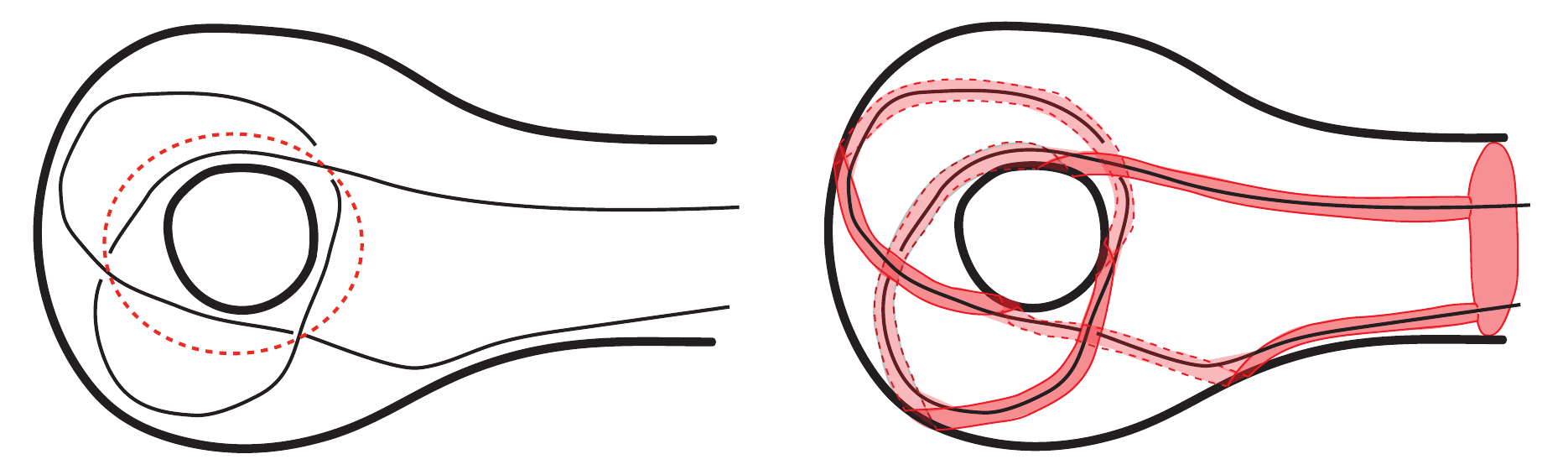}
    \caption{An essential annulus is present when the red curve does not satisfy condition (iv). 
    }
    \label{annulusharingcrossings}
\end{figure}

The conditions (i)-(iv) of Theorem \ref{thm-thickened} are necessary. Indeed, if (i) does not hold, then there is an essential twice-punctured sphere. If (ii) does not hold then a region with genus greater than 0 produces an essential annulus by taking $\alpha \times I$ for any nontrivial non-boundary-parallel simple closed curve $\alpha$ in the region. A planar region with more than one boundary produces an essential disk as in Figure \ref{fig-nonexamples-thickened}(a). 
If (iii) does not hold, then there is an essential annulus as in  Figure \ref{fig-nonexamples-thickened}(b). If (iv) does not hold, then we will show in  Lemma \ref{lemma-failiv-exists-annulus} that there is an essential annulus, an example of which appears in Figure \ref{annulusharingcrossings}. 
So our main task will be to prove that these conditions imply tg-hyperbolicity, which we  do in Section 3.

Note that if $F$ is a disk, then $M$ is the complement of a link in a 3-ball. Capping off the spherical boundary with a 3-ball yields a link complement in $S^3$. Then Corollary 2 of \cite{menasco} characterizes when $M$ is hyperbolic. Similarly, if $F$ is closed and of genus at least 1, Theorem 1 of \cite{small18} characterizes when $M$ is tg-hyperbolic.

As an application of Theorem \ref{thm-thickened}, note that $F \times I$ is always a handlebody. Specifically, if $F$ is an orientable genus $g$ surface with $k$ boundary components, then $F \times I$ is a genus $2g + (k - 1)$ handlebody. Hence, if $L$ is a link in a handlebody and we can find a way to represent this handlebody as a thickened projection surface $F$ such that $L$ is cellular alternating on $F$, we can determine if $L$ is hyperbolic. In particular, there are examples of links in handlebodies which do not have a closed projection surface satisfying the hypotheses of Theorem 1.1 from \cite{hp17}, but do have a projection surface with boundary satisfying the hypotheses of Theorem \ref{thm-thickened}. We give an example in Section \ref{section-apps}.

Relatively few links in handlebodies with tg-hyperbolic complements were previously known. One was proved to be so in \cite{adamsnew}. There are a finite number given in \cite{FMP}. Each of \cite{Frigerio} and \cite{simplesmallknots} give an explicit infinite set, one for each genus. Theorem 1.1 of \cite{hp17} does generate many examples when each compressing disk on the boundary of the handlebody is crossed at least four times by an appropriate alternating link.  Theorem \ref{thm-thickened}, particularly when conjoined with a form of ``composition'', as described in Theorem 2.1 of \cite{complinks}, further increases the number of such known.  

\subsection{Generalizing to additional ambient manifolds} 

We now consider a compact orientable 3-manifold with boundary $Y$ that contains a properly embedded orientable surface $F$ with boundary that is both incompressible  and $\partial$-incompressible and that intersects all essential annuli and tori in $Y$. We show that if an appropriate link is removed from $F \times I \subset Y$, the link complement in $Y$ is tg-hyperbolic.

First, we state a similarly easy way to check whether or not $L$ is prime in $Y$. Let $\pi(L)$ be a projection of a link $L$ to $F$.

\begin{proposition} \label{prop-weakly-prime-circ}
Let $F$ be an orientable incompressible $\partial$-incompressible surface with nonempty boundary properly embedded in a compact orientable irreducible $\partial$-irreducible 3-manifold $Y$, and let $L$ be a link in $Y$ with a connected, reduced, cellular
alternating 
projection to $F$. Then $L$ is prime in $Y$ if and only if $\pi(L)$ is weakly prime on $F$. 
\end{proposition}

Then we have the following theorem:

\begin{theorem} \label{surfaceinmanifold}
Let $Y$ be an orientable irreducible $\partial$-irreducible 3-manifold with $\partial Y \neq \emptyset$. Let $F$ be a properly embedded  orientable incompressible and $\partial$-incompressible connected surface with boundary in $Y$. Suppose all essential tori and annuli that exist in $Y$  intersect $F$. Let $L$ be a link in a regular neighborhood $N = F \times I$ of $F$ that has a connected reduced alternating projection $\pi(L)$ to $F$ with at least one crossing and that satisfies the following conditions:

\begin{enumerate} [label = (\roman*)]
    \item $\pi(L)$ is weakly prime on $F$; 
    \item the interior of every complementary region of $F  \setminus \pi(L)$ is either an open disk or an open annulus;

    \end{enumerate}

Then $Y \setminus N(L)$ is tg-hyperbolic. 
    
 \end{theorem}

 Compare to Theorem \ref{thm-thickened}; in particular, here, the previous condition (iii) that adjacent regions cannot both be annuli and the previous condition (iv) are no longer necessary. In either of these cases, the annulus that is generated when the condition fails does not extend to an annulus in $Y$.  Also, note again that condition (ii) implies the link is cellular alternating. 
 The requirement that $F$ be connected is not essential, as the theorem can be repeated for additional surface components and tg-hyperbolicity is maintained.

 Examples for $Y$ include any finite-volume hyperbolic 3-manifold with cusps and/or totally geodesic boundary of genus at least 2. In these cases, there are no essential annuli or tori to worry about and  there is always an incompressible $\partial$-incompressible surface that can play the role of $F$ (see for instance Lemma 9.4.6 of \cite{martelli}). A simple example for $F$ would be a minimal genus Seifert surface in a hyperbolic knot complement $Y$ in $S^3$. 

\begin{figure}[htbp]
    \centering
    \includegraphics[scale = 0.64]{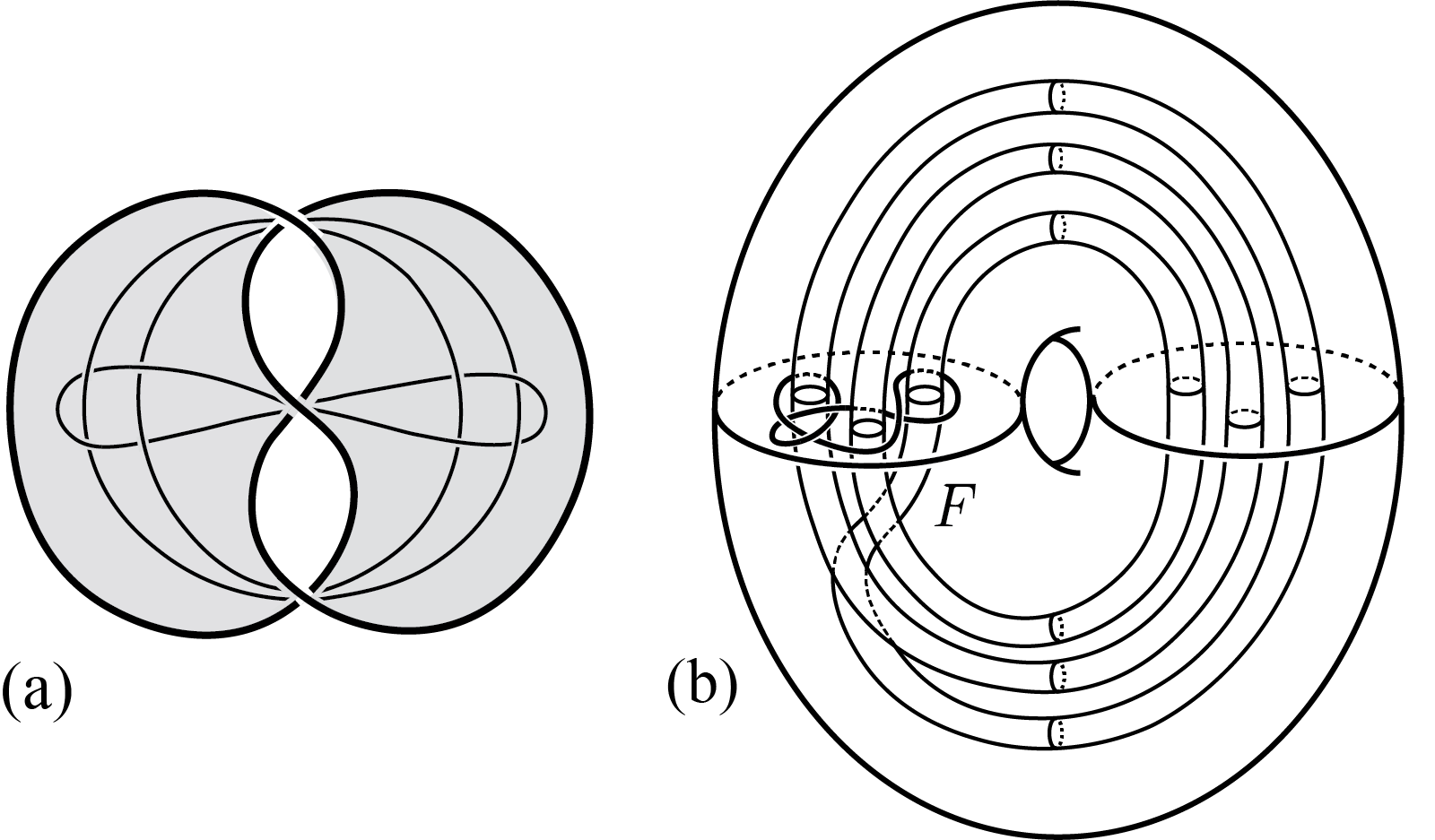}
    \caption{Examples of manifolds that are tg-hyperbolic by Theorem \ref{surfaceinmanifold}.}
    \label{Theorem1.8examples}
\end{figure}

 In Figure \ref{Theorem1.8examples}(a), we see an example where $Y$ is the complement of the trefoil knot. There is an essential annulus, however, it intersects the shaded Seifert surface playing the role of $F$. Hence the link complement shown is hyperbolic by Theorem \ref{surfaceinmanifold}.

Examples also include any surface bundle $F \Tilde{\times} S^1$ over a compact  orientable  surface $F$  with nonempty boundary other than a disk or annulus. We pick the incompressible $\partial$-incompressible surface to be a fiber. Such a manifold does contain essential annuli and/or tori but they all intersect $F$. See Figure \ref{Theorem1.8examples}(b) for an example.

\subsection{Organization and further directions}
In Section \ref{section-int-curves}, we introduce the notion of \textit{bubbles}, first defined by Menasco in \cite{menasco}. We also adapt various lemmas from \cite{small18} concerning intersection curves. In Section \ref{section-pf-thm-thickened}, we prove Theorem \ref{thm-thickened} as well as Proposition \ref{prop-weakly-prime}, and in Section \ref{section-pf-thm-circle-bundle}, we prove Theorem \ref{surfaceinmanifold} and Proposition \ref{prop-weakly-prime-circ}. 
Finally, in Section \ref{section-apps}, we discuss some applications of our results. One motivation for Theorem \ref{thm-thickened} is that it gives a large class of hyperbolic links in handlebodies, which we may naturally view as thickened surfaces-with-boundary. Besides being interesting objects in their own right, they appear naturally in the study of \textit{knotoidal graphs}, as defined in \cite{generalizedknotoids}. In that paper, the authors construct a map $\phi_\Sigma^D$ from the set of knotoidal graphs to the set of spatial graphs in 3-manifolds. Hence, there is a well-defined notion of hyperbolicity and hyperbolic volume for knotoidal graphs which may be used to distinguish them. 

In particular, \textit{staked links}, as defined in \cite{generalizedknotoids}, are obtained by adding vertices, called isolated poles, to the complementary regions of a link projection. As in the case of endpoints of knotoids, we do not allow strands of the projection to pass over or under these poles.   A subclass of knotoidal graphs, these staked links are mapped by $\phi_\Sigma^D$ into the set of links in handlebodies. Then it is interesting to determine which staked links are hyperbolic and compute their volumes. Theorem \ref{thm-thickened} gives the answer to the first question in the case of alternating staked links. Furthermore, we can show that certain staked links which are ``close to alternating" in some sense are also hyperbolic, using results from \cite{complinks}. Theorem \ref{thm-thickened} is used in \cite{generalizedknotoids} to prove that every link in $S^3$ can be staked to be hyperbolic.  

In future work, it would be interesting to obtain volume bounds for alternating links in handlebodies (and such results would immediately apply to alternating staked links). Also, while we work only with orientable 3-manifolds, we suspect similar results hold when $F$ is a nonorientable surface with boundary, or when we allow for orientation-reversing self-homeomorphisms of $F$ in the construction of $F \Tilde{\times} S^1$. 


\subsection*{Acknowledgements} The research was supported by Williams College and NSF Grant DMS-1947438 supporting the SMALL Undergraduate Research Project. We are grateful to Alexandra Bonat, Maya Chande, Maxwell Jiang, Zachary Romrell, Daniel Santiago, Benjamin Shapiro and Dora Woodruff, who are the other members of the knot theory group of the 2022 SMALL REU program at Williams College, for many helpful discussions and suggestions.

\label{section-intro}

\section{Bubbles and Intersection Curves} 
Let $F$ be a compact, connected, orientable surface with boundary, and let $L \subset F \times I$ be a link with connected, cellular alternating projection $\pi(L) \subset F \times \{1/2\} = F$. By Thurston's Hyperbolization Theorem, proving that $M$ has no essential spheres, tori, disks, or annuli is sufficient to conclude that $M$ is tg-hyperbolic.

Throughout this section, $\Sigma$ is a properly embedded essential surface in $M = (F \times I) \setminus N(L)$ with boundary on $\partial(F \times I)$. Note that when $\Sigma$ is a disk, we may always assume it has boundary on $\partial(F \times I)$. Indeed, suppose $\Sigma$ has boundary on $\partial N(K)$, where $K$ is a component of $L$. Then $\partial (N(K) \cup N(D))$ is an essential sphere. Hence, if we can eliminate essential spheres, then we have eliminated such disks. 

\begin{figure}[htbp]
    \centering
    \includegraphics[scale = 0.6]{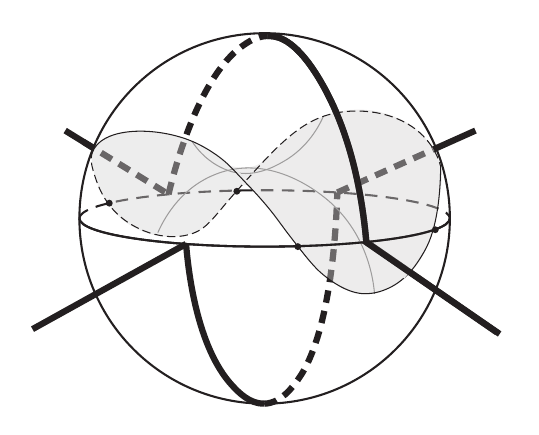}
    \caption{A surface $\Sigma$ intersecting a bubble in a saddle disk.
    }
    \label{fig-bubble-saddle}
\end{figure}

As in \cite{menasco}, arrange $L$ to lay in $F = F \times \{1/2\}$ away from the crossings, and at each crossing place a 3-ball $B$ which we call a \textit{bubble}. Arrange the over and understrands so that they lie in the upper hemisphere $\partial B_+$ and lower hemisphere $\partial B_-$ of the bubble, respectively. We may isotope $\Sigma$ to intersect the bubbles in saddle-shaped disks, by first isotoping $\Sigma$ to intersect the vertical axis of each bubble transversely, then pushing $\Sigma$ radially outward from the axis. See Figure \ref{fig-bubble-saddle}. 

Let $F_+$ (resp. $F_-$) be the surface obtained from $F \times \{1/2\}$ by removing each equatorial disk where $F \times \{1/2\}$ intersects a bubble $B$ and replacing it with the upper hemisphere $\partial B_+$ (resp. lower hemisphere $\partial B_-$). The desired contradictions come from analyzing the intersection curves between $\Sigma$ and the $F_\pm$, which may be closed or properly embedded arcs since $\partial \Sigma \subset \partial (F \times I)$ and therefore, can be perturbed to intersect the $F_\pm$ transversely. 


In the remainder of this section, we state and prove various lemmas for $\Sigma \cap F_+$, noting that all results and arguments apply to to $\Sigma \cap F_-$ as well. 

\begin{lemma} \label{lemma-at-least-one-int-curve}
There is at least one intersection curve in $\Sigma \cap F_+$. 
\end{lemma}

\begin{proof}
Suppose for contradiction that $\Sigma \cap F_+ = \emptyset$. Without loss of generality, $\Sigma$ can be isotoped to lie in $F \times (1/2, 1]$, a handlebody, and if $\Sigma$ has boundary, then $\partial \Sigma$ lies in $(F \times \{1\}) \cup (\partial F \times (1/2, 1])$. If $\Sigma$ is a sphere or torus, then it is compressible, and if $\Sigma$ is a disk or annulus, then it is $\partial$-parallel. 
\end{proof}

We would like to simplify $\Sigma \cap F_+$ as much as possible. Assign to each embedding of $\Sigma$ an ordered pair $(s, i)$, where $s$ is the number of saddle disks in the intersection between $\Sigma$ and the bubbles and $i$ is the number of intersection curves in $\Sigma \cap F_+$. For the remainder of the section, we may assume that our choice of an embedding of $\Sigma$ minimizes $(s, i)$ under lexicographical ordering. 


To this end, we can show that $\Sigma \cap F_+$ cannot contain any intersection curves which are \textit{trivial} on both $\Sigma$ and $F_+$. 

\begin{definition}
We say a simple closed curve on a surface is \textit{trivial} if it bounds a disk in the surface. We say a properly embedded arc on a surface is \textit{trivial} if it cuts a disk from the surface. 
\end{definition}

To eliminate curves trivial on both $\Sigma$ and $F_+$, we define the notion of \textit{meridional (in)compressibility}, first introduced in \cite{menasco}. Eventually, in Section \ref{section-pf-thm-thickened}, we show that essential surfaces in $M$ cannot be meridionally incompressible nor meridionally compressible, thus eliminating them.


\begin{definition}
Let $Y$ be a compact 3-manifold containing a link $L$, and let $\Sigma$ be a properly embedded surface in $Y \setminus N(L)$. We say $\Sigma$ is \textit{meridionally incompressible} if for every disk $D$ in $Y$ such that $D \cap \Sigma = \partial D$ and $D$ is punctured exactly once by $L$,  there is another disk $D'$ in $\Sigma \cup N(L)$ such that $\partial D' = \partial D$ and $D'$ is punctured by $L$ exactly once. Otherwise, we say $\Sigma$ is \textit{meridionally compressible} and we call $D$ a \textit{meridional compression disk}. We refer to surgery on $\Sigma$ along $D$ as a \textit{meridional compression}. 
\end{definition}

Throughout this section we take $Y = F \times I$. We remark that if $\Sigma$ is essential in $M$ and has a meridional compression disk $D$, then the (not necessarily connected) surface $\Sigma'$ resulting from meridionally compressing $\Sigma$ along $D$ is also essential. 

The following lemma tells us that when $\Sigma$ is meridionally incompressible, all of the closed intersection curves are nontrivial in some sense. 

\begin{lemma} \label{lemma-removing-triv-curve}
Suppose $\Sigma$ is meridionally incompressible and $(s, i)$ is minimized. Then no closed intersection curve in $\Sigma \cap F_+$ can be trivial in both $\Sigma$ and  $F_+$. 
\end{lemma}

\begin{figure}[htbp]
    \centering
    \includegraphics[scale = 0.6]{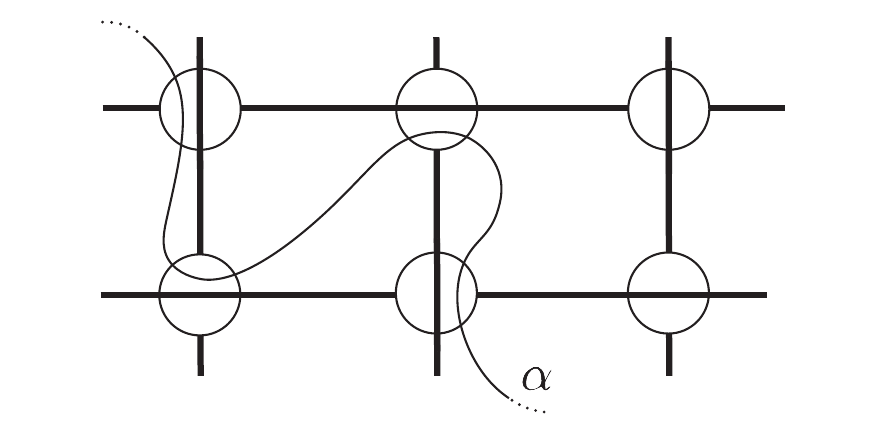}
    \caption{The alternating property: $\alpha$ alternates sides as it encounters each bubble. }
    \label{fig-alternating-property}
\end{figure}

In order to prove this, we first prove several other lemmas. Note that the upper hemisphere of each bubble $B$ is separated by the overstrand into two sides. We observe that because $\pi(L)$ is an alternating diagram, any intersection curve in $\Sigma \cap F_+$ must alternate between entering bubbles on the left side of the overstrand and entering on the right side of the overstrand. See Figure \ref{fig-alternating-property} for a local picture. The alternating property places strong restrictions on the appearance of the intersection curves in $F_+$, as the next two lemmas demonstrate. 

\begin{lemma} \label{lemma-same-side-bubble}
A closed intersection curve $\alpha$ in $\Sigma \cap F_+$ which is trivial in $F_+$ cannot intersect a bubble twice on the same side. 
\end{lemma}

\begin{figure}[htbp]
    \centering
    \includegraphics[scale = 1.1]{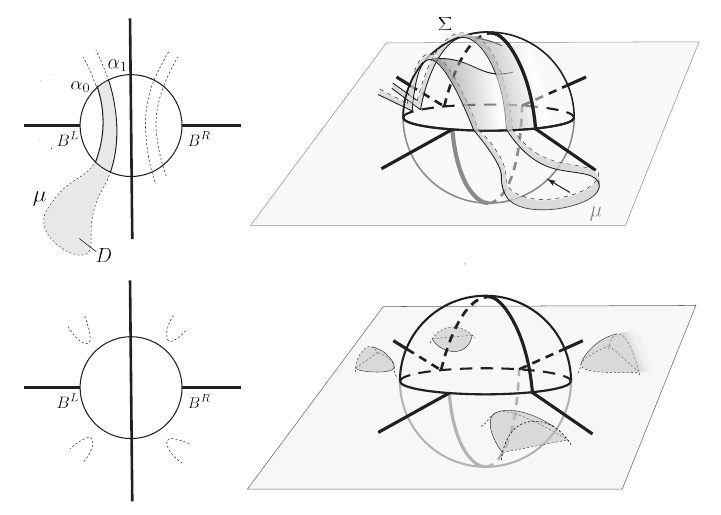}
    \caption{Isotoping $\Sigma$ to remove two saddles. We visualize $\Sigma$ as folding back over itself to intersect $B^L$ twice. Then imagine using your finger to push the fold into the bubble and downwards past $F$, smoothing it out and removing $\alpha_0$ and $\alpha_1$. For simplicity, we only draw $\Sigma$ where it lies above $F$.  
    }
    \label{fig-finger-saddle-move}
\end{figure}

\begin{proof}
Let $B$ be a bubble and for convenience, let $B^L$ and $B^R$ denote the two halves of $B$ obtained by slicing along a vertical plane containing the overstrand. Without loss of generality, suppose $\alpha$ meets $B^L$ in more than one arc, and let $D$ be the disk in $F_+$ bounded by $\alpha$. Let $\{\alpha_i\}$ be the set of arcs in $\alpha \cap B^L$ and observe that each corresponds to a distinct saddle disk in $\Sigma \cap B$. Furthermore, we may assume that there exists a pair of arcs, $\alpha_0$ and $\alpha_1$, which are adjacent on $B^L$. Indeed, if there is another intersection curve $\alpha'$ in $\Sigma \cap F_+$ which meets $B^L$ in between $\alpha_0$ and $\alpha_1$, then it does so in at least two arcs between $\alpha_0$ and $\alpha_1$ because $\alpha'$ must be contained in $D$. We can continue finding these nested loops until we are able to choose an adjacent pair of arcs belonging to the same intersection curve. 

As in \cite{adams}, let $\mu$ be an arc running along $\alpha$ in $\Sigma$ connecting $\alpha_0$ and $\alpha_1$. Then we may isotope $\Sigma$ to remove the two saddles corresponding to $\alpha_0$ and $\alpha_1$ by pushing a regular neighborhood of $\mu$ across the disk $D$, through $B$, and downwards past $F_+$. See Figure \ref{fig-finger-saddle-move}. This reduces $s$ by $2$, contradicting minimality of $(s, i)$. 
\end{proof}

\begin{lemma} \label{lemma-both-sides-bubble}
Suppose $\Sigma$ is meridionally incompressible and suppose $\alpha$ is a closed intersection curve in $\Sigma \cap F_+$ which is trivial in $F_+$ and bounds a disk containing the overstrand of some bubble $B$. Then $\alpha$ cannot intersect $B$ on both sides.
\end{lemma}


\begin{proof}
Suppose for contradiction that a curve $\alpha$ satisfying the hypotheses intersects $B$ on both sides in arcs $\alpha_0$ and $\alpha_1$. Let $D$ be a disk in $F_+$ bounded by $\alpha$. 

First, observe that $\alpha_0$ and $\alpha_1$ must be connected by an arc $\gamma$ as in Figure \ref{fig-both-sides}(a).  Otherwise, suppose $\gamma$ appears as in Figure \ref{fig-both-sides}(b). Then because $\alpha$ bounds a disk, the endpoint of $\alpha$ cannot escape the enclosed region along a handle and hence must exit the enclosed region by passing over $B$ again. However, whichever side of $B$ it passes over, this contradicts Lemma \ref{lemma-same-side-bubble}.

\begin{figure}[htbp]
    \centering
    \includegraphics[scale = .55]{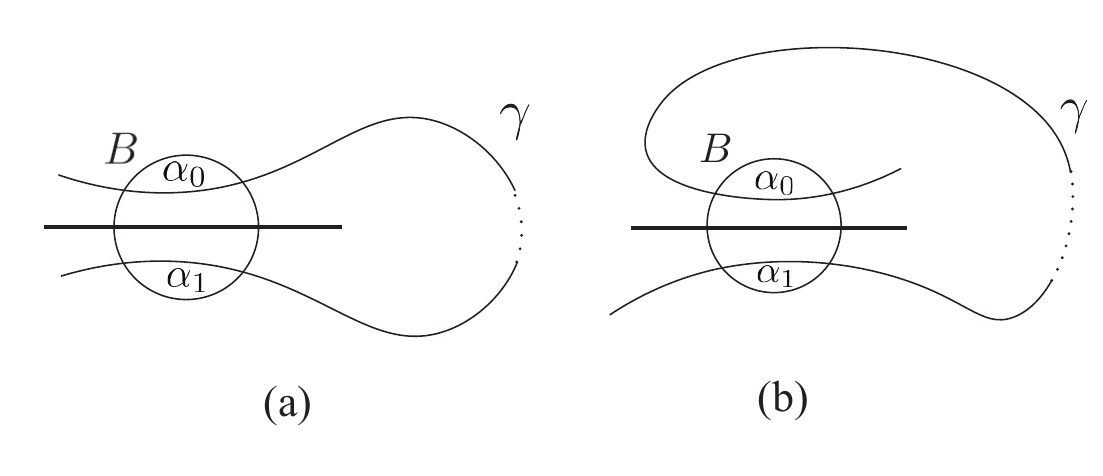}
    \caption{The two possibilities when an intersection curve intersects the opposite sides of a bubble. We rule out 8(b) by contradiction with Lemma \ref{lemma-same-side-bubble}. 
    }
    \label{fig-both-sides}
\end{figure}

So now we can assume we are in a situation as depicted in Figure \ref{fig-both-sides}(a). Of the intersection curves which are trivial and intersect $B$ on both sides, we may choose $\alpha$ to intersect closest to the overstrand of $B$ on one of its sides. We claim that the two arcs of $\alpha$ which intersect $B$ closest to the overstrand on either side belong to the same saddle disk of $\Sigma \cap B$. Indeed, suppose $\alpha'$ is another intersection curve such that it meets $B$ between an arc of $\alpha \cap B$ and the overstrand. Since $\alpha$ bounds a disk, $\alpha'$ must intersect $B$ at least twice on one side of $B$. Since $D$ contains the overstrand of $B$, $\alpha'$ is contained in $D$ and hence, is trivial. This contradicts Lemma \ref{lemma-same-side-bubble}. 

Let $\sigma \subset \Sigma \cap B$ be the saddle disk  which contains $\alpha_0$ and $\alpha_1$ in its boundary. The remainder of the proof more or less follows the proof of Lemma 7 of \cite{small18} or Lemma 1 of \cite{menasco}. Let $\mu$ be an arc in $\Sigma$ running parallel to $\alpha$ which connects $\alpha_0$ and $\alpha_1$. Then there is an arc $\gamma$ in $\sigma$ such that  $\mu \cup \gamma$ is a circle in $\Sigma$ which bounds a  disk punctured once by $L$ in $M$. See Figure \ref{fig-both-sides-comp-disk}. But this yields a meridional compression disk, a contradiction. 

\begin{figure} [htbp]
    \centering
    \includegraphics[scale=.8]{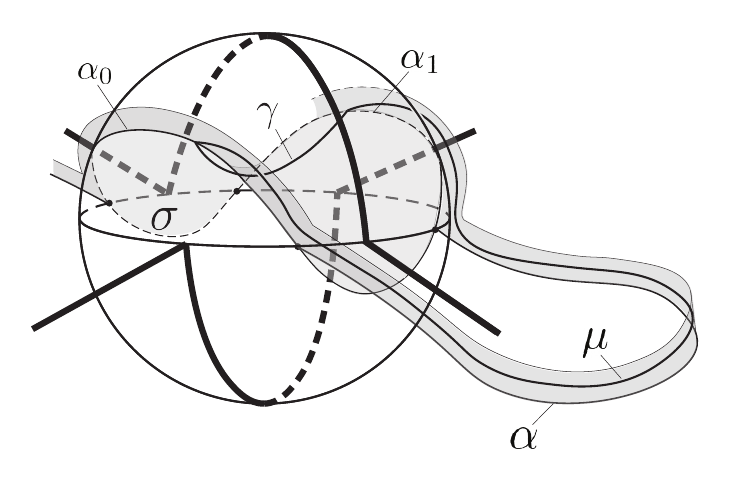}
    \caption{The curve $\mu \cup \gamma$ bounds a meridional compression disk for $\Sigma$.
    }
    \label{fig-both-sides-comp-disk}
\end{figure}
\end{proof}

To prove Lemma \ref{lemma-removing-triv-curve} and subsequent lemmas, the notion of innermost (closed) curves and outermost arcs is very useful. 

\begin{definition}
Let $\alpha$ be a closed intersection curve in $\Sigma \cap F_+$. If $\alpha$ is trivial on $\Sigma$ (resp. trivial on $F_+$), we say $\alpha$ is \textit{innermost} on $\Sigma$ (resp. on $F_+$) if $\alpha$ bounds a disk $D$ in $\Sigma$ (resp. $F_+$) which does not contain any intersection curves of $\Sigma \cap F_+$ in its interior. Similarly, suppose $\alpha$ is an intersection arc in $\Sigma \cap F_+$. If $\alpha$ is trivial on $\Sigma$ (resp. $F_+$), we say $\alpha$ is \textit{outermost} on $\Sigma$ (resp. on $F_+$) if $\alpha$, together with an arc of $\partial \Sigma$ (resp. $\partial F_+$), bounds a disk $D$ in $\Sigma$ (resp. $F_+$) which does not contain any intersection curves of $\Sigma \cap F_+$ in its interior.
\end{definition}

\begin{proof}[Proof of Lemma \ref{lemma-removing-triv-curve}]
Let $D \subset F_+$ and $D' \subset \Sigma$ be disks bounded by $\alpha$, and of the (closed) intersection curves contained in $D$, let $\beta$ be an innermost such curve on $F_+$. Then $\beta$ bounds a disk $E \subset D \subset F_+$ and by Lemmas \ref{lemma-same-side-bubble} and \ref{lemma-both-sides-bubble}, $\beta$ cannot intersect any bubbles. By iterating this argument, we find that none of the intersection curves contained in $D$ can intersect bubbles and hence, neither can $\alpha$. Then $D \cup D'$ is a 2-sphere which is not punctured by the link. By irreducibility of $F \times I$, $D \cup D'$ bounds a 3-ball and we may isotope $D'$ through the 3-ball, pushing it slightly past $F_+$ to remove $\alpha$ and any intersection curves contained in $D$. This contradicts minimality of $(s, i)$. 
\end{proof}

We conclude this section by using Lemma \ref{lemma-removing-triv-curve} to prove several more lemmas restricting the appearance of the intersection curves when $\Sigma$ is meridionally incompressible.

\begin{lemma} \label{lemma-triv-int-curve-iff}
Suppose $\Sigma$ is meridionally incompressible, and let $\alpha \subset \Sigma \cap F_+$ be an intersection curve. Then $\alpha$ is trivial in $\Sigma$ if and only if one of the following is true: 
\begin{enumerate}
    \item $\alpha$ is trivial in $F_+$; or
    \item $\alpha$ is an arc bounding a $\partial$-compression disk for $F_+$ in $M$.
\end{enumerate}
\end{lemma}

\begin{proof} 
Suppose $\alpha$ is nontrivial in $\Sigma$ and trivial in $F_+$. When $\Sigma$ is a sphere or disk, this is clearly impossible. Of all intersection curves which are nontrivial in $\Sigma$ and trivial in $F_+$, choose $\alpha$ to be an innermost closed curve or an outermost arc on $F_+$. Let $D$ be a disk in $F_+$ bounded by $\alpha$ (and possibly an arc of $\partial F_+$ if $\alpha$ is an arc). Then any intersection curve in $D$ is trivial in both $\Sigma$ and $F_+$, so by Lemma \ref{lemma-removing-triv-curve}, we may assume that all of these curves are arcs. 


In particular, if $\Sigma$ is a torus, or $\Sigma$ is an annulus and $\alpha$ is a closed curve, then all intersection curves contained in $D$ are eliminated, and $\alpha$ bounds a compression disk for $\Sigma$, a contradiction. If $\Sigma$ is an annulus and $\alpha$ is an arc, then all closed intersection curves in $D$ are eliminated. If all intersection curves in $D$ are closed, then $\alpha$ bounds a compression disk for $\Sigma$. Hence, there is at least one intersection arc contained in $D$. But then the outermost such arc bounds a $\partial$-compression disk for $\Sigma$, contradicting essentiality. 

Conversely, suppose $\alpha$ is trivial in $\Sigma$ and nontrivial in $F_+$. Fill in $N(L)$ to work in the handlebody $F \times I$. We show that $\Sigma$ cannot intersect $F_+$ in such an $\alpha$ in $F \times I$, much less in $(F \times I) \setminus N(L)$. Of all the curves which are trivial in $\Sigma$ and nontrivial in $F_+$, choose $\alpha$ to be an innermost closed curve or outermost arc on $F_+$, and let $D'$ be a disk in $\Sigma$ bounded by $\alpha$ (and possibly an arc of $\partial \Sigma$ if $\alpha$ is an arc). Every intersection curve contained in $D'$ is trivial in $\Sigma$ and in $F_+$. Of these curves, suppose at least one is closed and choose $\beta$ to be an innermost such curve on $\Sigma$. Let $E$ and $E'$ be disks bound by $\beta$ in $F_+$ and $\Sigma$ respectively. Push the interior of $E$ slightly off $F_+$ so $E \cup E'$ is a 2-sphere in $M$. By irreducibility of $F \times I$, $E \cup E'$ bounds a 3-ball. Isotope $\Sigma$ through this 3-ball to remove $\beta$ as before, and by iterating this process, we can remove all closed intersection curves contained in $D$. If $\alpha$ is a closed curve, then all intersection curves in $D$ are eliminated, and $\alpha$ bounds a compression disk for $F_+$ in $F \times I$, a contradiction. If $\alpha$ is an arc, it is either trivial in $F_+$ or $D'$ is a $\partial$-compression disk for $F_+$ in $F \times I$, hence in $M$.
\end{proof}

Henceforth, we say an intersection curve is trivial if it is trivial on $\Sigma$ and $F_+$. When $\Sigma$ is meridionally incompressible we may assume such curves are arcs by Lemmas \ref{lemma-removing-triv-curve} and \ref{lemma-triv-int-curve-iff}.

\begin{lemma} \label{lemma-int-curve-one-bubble}
Suppose $\Sigma$ is meridionally incompressible, and let $\alpha \subset \Sigma \cap F_+$ be an intersection arc. 
Then $\alpha$ intersects at least one bubble. 
\end{lemma}

\begin{proof}

Suppose for contradiction that $\alpha$ does not intersect any bubble. Then the arc $\alpha$ is properly embedded in a complementary region $R$ of $F_+ \setminus \pi(L)$ which is homeomorphic to an annulus, and furthermore, both endpoints of $\alpha$ lie on the same boundary component of $R$. In particular, $\alpha$ is trivial in $F_+$ and $\Sigma$, using Lemma \ref{lemma-triv-int-curve-iff}. 

Let $D$ be a disk in $R \subset F_+$ bounded by $\alpha$ and an arc $\beta$ of $\partial F_+$, and let $D'$ be a disk in $\Sigma$ bounded by $\alpha$ an arc $\gamma$ of $\partial \Sigma \subset \partial (F \times I)$. Then $D \cup D'$ is a disk with boundary $\beta \cup \gamma$ on $\partial (F \times I)$. Furthermore, since $\alpha$ is outermost on $\Sigma$, $\gamma$ does not intersect $\partial F_+$ away from its endpoints. Hence, without loss of generality, $\beta \cup \gamma$ lies in $(F \times \{1\}) \cup (\partial F \times [1/2, 1])$. But $D \cup D'$ cannot be a compression disk for $F \times [1/2, 1]$, so $\beta \cup \gamma$ must bound a disk $E$ in $(F \times \{1\}) \cup (\partial F \times [1/2, 1])$. Then $D \cup D' \cup E$ is a 2-sphere in $F \times [1/2, 1]$ bounding a 3-ball. Isotope $D'$ to $D$ through the 3-ball and slightly past $F_+$, removing $\alpha$ and any other intersection curves contained in $D$. This reduces $i$ without affecting $s$, a contradiction.

\end{proof}

\label{section-int-curves}

\section{Proof of Theorem \ref{thm-thickened}} 
For now, we assume the statement Proposition \ref{prop-weakly-prime} so we may regard $L$ as prime in $F \times I$. We prove that $M$ has no essential spheres, tori, disks, or annuli in that order. To use the lemmas from the previous section, we assume that $\partial \Sigma \subset \partial(F \times I)$: later on we eliminate essential annuli with at least one boundary component on $\partial N(L)$ using different methods.

We will be explicit about which of the conditions we use from Theorem \ref{thm-thickened} so that we can also use the appropriate lemmas for the proof of Theorem \ref{surfaceinmanifold}.

\begin{lemma} \label{lemma-torus-annulus-triv-curve}
Let $\Sigma$ be a meridionally incompressible 
essential torus or an essential annulus with both boundary components on $\partial(F \times I)$ such that $L$ satisfies the hypotheses of Theorem \ref{thm-thickened} with the possible exception of conditions (iii) and (iv). Then at least one intersection curve in either $\Sigma \cap F_+$ or $\Sigma \cap F_-$ is trivial on $\Sigma$. 
\end{lemma}

\begin{proof} 
Consider $\Sigma$ with all of the intersection curves in $\Sigma \cap F_+$ and $\Sigma \cap F_-$ projected onto it along with saddles corresponding to quadrilaterals. After shrinking down the saddles to vertices, we obtain a (not necessarily connected) 4-valent graph together with some circles without vertices.

These circles correspond to nontrivial closed intersection curves which do not meet any bubbles and lie in an annular region of $F_+ \setminus \pi(L)$. Denote the graph together with the circles without vertices by $\Gamma$; we call $\Gamma$ the \textit{intersection graph} of $\Sigma$ (see Figure \ref{fig-int-graph}). Each complementary region of $\Sigma \setminus \Gamma$ lies entirely in one of the components of $M \setminus F$, and the boundary of a region is precisely an intersection curve in either $\Sigma \cap F_+$ or $\Sigma \cap F_-$, depending on which side of $F$ the region lies in. Hence, our goal is to show that at least one region of $\Sigma \setminus \Gamma$ is a disk.

\begin{figure}[htbp]
    \centering
    \includegraphics[scale=.45]{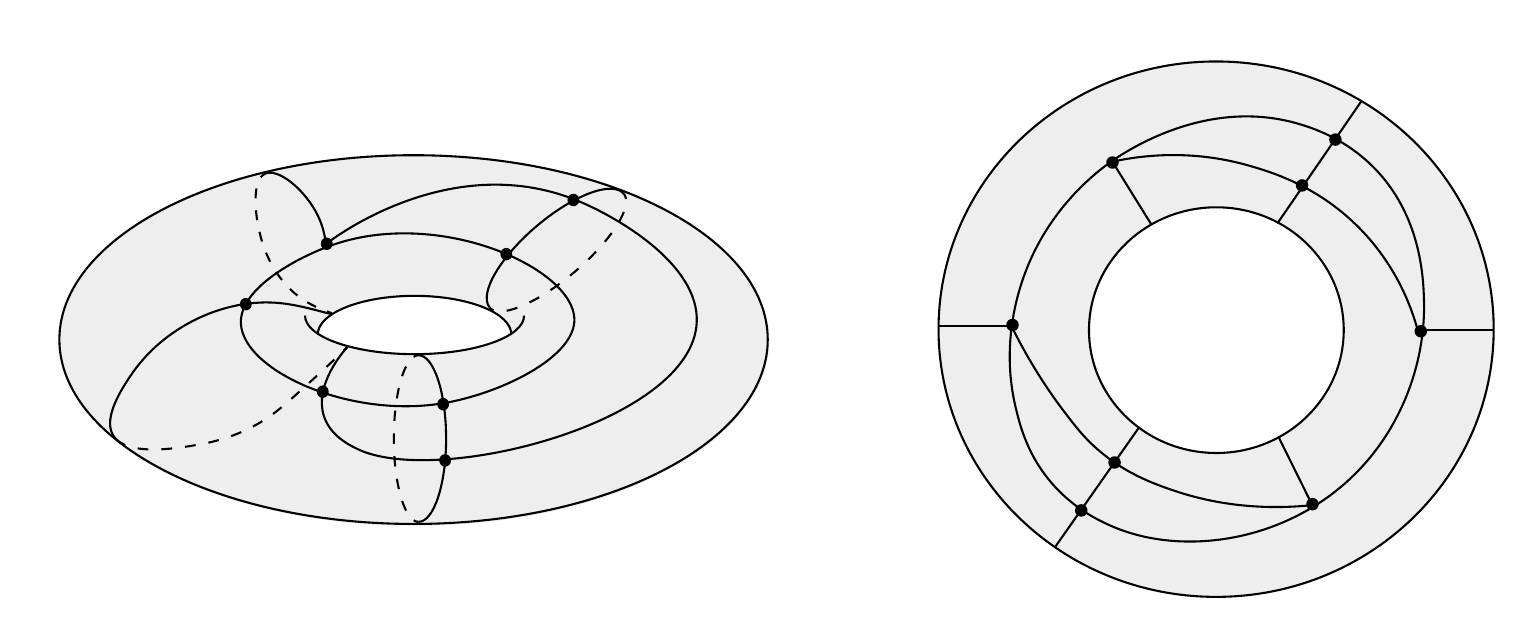}
    \caption{Possible intersection graphs for $\Sigma$.}
    \label{fig-int-graph}
\end{figure}

First, suppose there are no circles without vertices in $\Gamma$,  so $\Gamma$ is just a 4-valent properly embedded graph on $\Sigma$. Let $V$, $E$, and $F$ be the number of saddles, edges, and regions of $\Sigma \setminus F$ respectively, including endpoints of intersection arcs on the boundary of $\Sigma$ as vertices. The Euler characteristic $\chi(\Sigma)$ is 0 when $\Sigma$ is a torus or annulus. From Lemmas \ref{lemma-at-least-one-int-curve} and \ref{lemma-int-curve-one-bubble}, we know there are  vertices from saddles. 
If no regions are disks, then each has nonpositive Euler characteristic contribution. 

On the torus, $E = 2V$, and we then have 
$$
\chi(\Sigma) \leq V - E = V - 2V = -V < 0. 
$$
Hence, when there are no circles without vertices, there is a disk region of $\Sigma \setminus \Gamma$, and its boundary will be an intersection curve in $\Sigma \cap F_+$ or $\Sigma \cap F_-$ which is trivial on $\Sigma$. 

On the annulus, let $V_I$ be the number of interior vertices corresponding to saddles and $V_B$ be the number of boundary vertices. Then 

$$\chi(\Sigma) \leq V_I + V_B  - \frac{4 V_I + 3 V_B}{2} = -V_I - \frac{V_B}{2} < 0.$$

Again, this forces there to be a region with positive Euler characteristic, meaning there is a disk region.


Now suppose that $\Gamma$ has at least one circle without vertices and at least one vertex. Note that the circles without vertices must be parallel to one another. Hence, the circles cut $\Sigma$ into annuli, and at least one annulus $A$ contains a nonempty 4-valent graph $\Gamma'$. Then the same Euler characteristic argument shows that at least one region of $A \setminus \Gamma'$ is a disk. 

It remains to show that $\Sigma$ must meet the bubbles in at least one saddle. 
Suppose $\Sigma$ does not meet any bubbles for contradiction, and consider the projection of the intersection curves in $\Sigma \cap F$ onto $F$. Then the intersection curves must be contained in one annular region of $F \setminus \pi(L)$ and encircle the same boundary component of $F$. First consider the case where $\Sigma$ is a torus. Then there are an even number of curves since $F$ is separating in $M$ (in particular, there are at least two). Take a pair of intersection curves $C_1$ and $C_2$ which are adjacent in $\Sigma$ (note they might not be adjacent in $F_+$). They bound an annulus $A$ in $\Sigma$ and an annulus $A'$ in $F_+$. 

After pushing $A'$ slightly off $F_+$, we find a torus $T = A \cup A'$ contained in one component of $M \setminus F$, both components of which are homeomorphic to a handlebody $F \times I$. Then observe that $T$ must have a compressing disk, as $\pi_1(T) \cong \mathbb{Z}^2$ cannot inject into $\pi_1(F \times I) \cong \mathbb{Z}^{\ast g}$. Compress $T$ along some compressing disk to obtain a 2-sphere, which bounds a 3-ball. Gluing back in the compressing disk yields a solid torus bounded by $T$. It cannot yield a knot exterior because of the fact $A'$ lies in the boundary of the handlebody. Note that the boundary curves of  the two annuli must intersect the boundary of the compressing disk once. So we can isotope $\Sigma$ through this solid torus to remove intersection curves $C_1$ and $C_2$, contradicting minimality of $(s, i)$. 

If $\Sigma$ is an annulus and there are an even number of intersection curves, the same argument goes through. Otherwise, we may assume $\Sigma$ intersects $F_+$ in a single intersection curve. Then we claim $\Sigma$ is $\partial$-parallel. 

Let $\gamma_0$ and $\gamma_1$ denote the boundary circles of $\Sigma$: since $\Sigma \cap F_+$ contains a single curve $\gamma_{\frac{1}{2}}$, we know that $\gamma_1$ must be in $(\partial F \times (\frac{1}{2}, 1]) \cup (F \times \{1\})$ and $\gamma_0$ must be in $(\partial F \times [0, \frac{1}{2})) \cup (F \times \{0\})$ (note that both of these spaces are homeomorphic to $F$). Viewed on $F$, the $\gamma_i$ for $i = 0, \frac{1}{2}, 1$ are all homotopic to the same component of $\partial F$. This means that $\partial \Sigma$ cuts an annular component $A$ from $\partial (F \times I)$. 
As above, $\Sigma$ must be parallel to $A$ making it $\partial$-parallel,  a contradiction to its being essential.

\end{proof}

Now we are ready to eliminate essential spheres and tori. 

\begin{proposition} \label{prop-no-essential-spheres-tori}
Let $M = (F \times I) \setminus N(L)$ as in Theorem \ref{thm-thickened} with the possible exception of conditions (iii) and (iv). Then $M$ has no essential spheres or tori. 
\end{proposition}

\begin{proof}
Suppose $\Sigma$ is meridionally compressible. If $\Sigma$ is a 2-sphere, then a meridional compression yields two 2-spheres in $F \times I$, each punctured by $L$ exactly once, which cannot occur. If $\Sigma$ is a torus, then a meridional compression yields a 2-sphere which is twice-punctured by the link. By Proposition \ref{prop-weakly-prime}, $L$ is prime in $F \times I$. Then $\Sigma$ must be $\partial$-parallel. 

Hence, $\Sigma$ must be meridionally incompressible. By Lemmas \ref{lemma-torus-annulus-triv-curve} and \ref{lemma-triv-int-curve-iff}, we may assume without loss of generality that $\Sigma \cap F_+$ contains some intersection curve $\alpha$ which is trivial in $\Sigma$. But by Lemma \ref{lemma-triv-int-curve-iff}, $\alpha$ is also trivial in $F_+$, and this contradicts Lemma \ref{lemma-removing-triv-curve}. Let $D$ be a disk in $F_+$ bounded by $\alpha$, and choose $\beta$ to be the innermost intersection curve contained in $D$. Let $D'$ be a disk in $F_+$ bounded by $\beta$. Since $\beta$ is trivial and closed, it intersects at least two bubbles (counted with multiplicity). Because of the alternating property, $D'$ contains the overstrand of some bubble $B$ intersected by $\beta$, and since $D'$ contains no other intersection curve, we conclude that $\beta$ must intersect $B$ on both sides (see Figure \ref{fig-no-essential-S2T2}). But this contradicts Lemma \ref{lemma-both-sides-bubble}.

\begin{figure}[htbp]
    \centering
    \includegraphics[scale = .75]{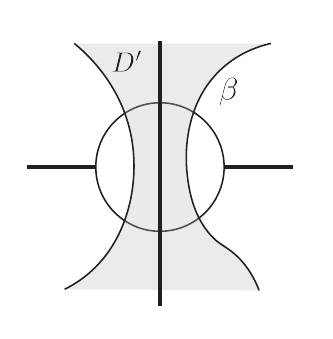}
    \caption{If $\Sigma$ is a meridionally incompressible essential sphere or torus, then there must be some curve $\beta$ bounding $D' \subset F_+$ appearing in this configuration at some bubble.}
    \label{fig-no-essential-S2T2}
\end{figure}
\end{proof}

As previously remarked, this also eliminates the possibility of an essential disk whose boundary lies on $N(K)$ for some component $K$ of $L$. 
To eliminate essential disks in general, we introduce the useful notion of \textit{forks} (which are similar to forks as defined in \cite{adams-forks} but here they only have two prongs). Consider an essential disk $\Sigma$ with the intersection curves in $\Sigma \cap F_+$ projected onto it. After shrinking down the vertices, we obtain a 4-valent intersection graph $\Gamma$ on $\Sigma$. 

\begin{definition} \label{defn-forks}
A \textit{fork} of $\Gamma$ is a vertex with at least two non-opposite edges ending on $\partial \Sigma$.  
\end{definition}

\begin{figure}[htbp]
    \centering
    \includegraphics[scale = .7]{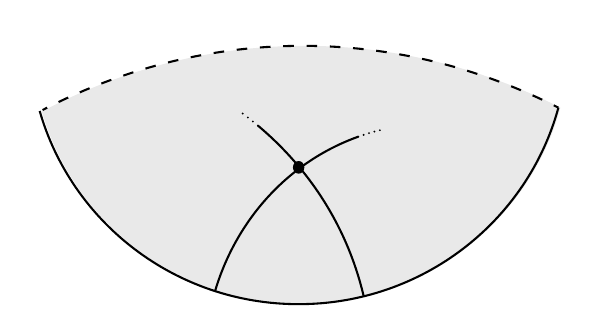}
    \caption{A portion of the intersection graph on $\Sigma$ corresponding to a fork.}
    \label{fig-fork-in-M}
\end{figure}

See Figure \ref{fig-fork-in-M}. Note that the endpoints of the two edges need not be adjacent on the boundary of $\Sigma$.



\begin{proposition} \label{prop-no-essential-disks}
Let $M = (F \times I) \setminus N(L)$ as in Theorem \ref{thm-thickened} with the possible exception of condition (iv). Then $M$ has no essential disks. 
\end{proposition}

\begin{proof}
First, note that suppose $\Sigma$ cannot be meridionally compressible as a meridional compression would generate a disk with boundary on a meridian of $N(L)$, that is, a 2-sphere once-punctured by the link, which cannot occur. 

Hence, $\Sigma$ must be meridionally incompressible, and by Lemma \ref{lemma-int-curve-one-bubble}, the intersection graph $\Gamma$ contains at least one vertex. As in \cite{adams-forks}, we can show there is at least one fork in $\Gamma$ as follows. Because there are no closed intersection curves using Lemmas \ref{lemma-removing-triv-curve} and \ref{lemma-triv-int-curve-iff}, every complementary region on the disk must intersect the boundary of the disk. If we discard all edges that touch the boundary to obtain a new graph $\Gamma'$ on $\Sigma$, we are left with a collection of trees. If a tree has two or more vertices, then it must have two or more leaves, each of which has three edges on $\Gamma$ that end on the boundary. So we have a fork. 

If the tree is only a single vertex, then there are four edges leaving it that end on the boundary of the disk and we again have a fork.

Let $\alpha$ and $\beta$ denote two adjacent edges of the fork. Note that $\alpha$ and $\beta$ cannot have endpoints on the same component of $\partial F$; otherwise, $\alpha \cup \beta$ together with an arc in $\partial F$ and an arc of the saddle bound a disk in $F_+$ which has exactly one intersection point with $\pi(L)$. Since $\alpha$ and $\beta$ don't meet other bubbles, they lie in adjacent complementary regions. But this implies that two adjacent regions are homeomorphically annuli, contradicting condition (iii) in the statement of Theorem \ref{thm-thickened}. 
\end{proof}

\begin{figure} [htbp]
    \centering
    \includegraphics[scale=.6]{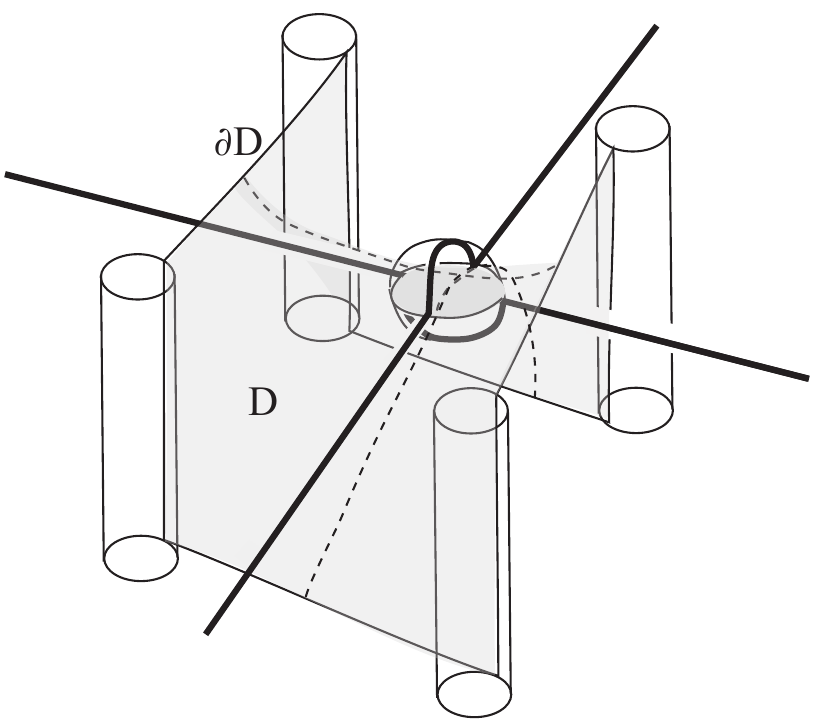}
    \caption{There is an essential disk if we allow all four complementary regions meeting at a crossing bubble to be annuli: it meets exactly one bubble in exactly one saddle disk.}
    \label{fig-essential-disk-wo-ciii}
\end{figure}

We remark that condition (iii) (or perhaps some weaker variation thereof) is necessary to show that no disks exist: for example, if all four complementary regions meeting at a crossing bubble are annuli, then there is an essential disk $D$ as pictured in Figure \ref{fig-essential-disk-wo-ciii}. 

It remains to show that there are no essential annuli. Suppose $\Sigma$ is an essential annulus. There are three cases to consider: 
\begin{enumerate}
    \item $\Sigma$ has both boundary components on $\partial (F \times I)$; 
    \item $\Sigma$ has one boundary component on $\partial (F \times I)$ and one boundary component on $\partial N(L)$; 
    \item $\Sigma$ has both boundary components on $\partial N(L)$. 
\end{enumerate}

\begin{lemma} \label{lemma-type2and3annuliimplytype1} Let $M = (F \times I) \setminus N(L)$ as in Theorem \ref{thm-thickened} with the possible exception of conditions (iii) and (iv). Then if $M$ has a essential annulus of type (2) or type (3), it has an essential torus or an essential annulus of type (1). 
\end{lemma}

\begin{proof} Let $A$ be an essential annulus of type (3), so both of its boundaries are on $\partial N(L)$. If there exists a single component $K$ such that $\partial A \subset \partial N(K)$, then let $T_1$ and $T_2$ be the two tori that form the boundary of $N(A \cup N(K))$. 

The boundaries of $A$ cannot be meridians on $\partial N(K)$, as that would contradict primeness. The annulus $A$ cuts $(F \times I) \setminus N(K)$ into two components, one of which contains $\partial (F \times I)$. Choose $T_1$ to be the torus that separates $K$ from $\partial (F \times I)$.   We prove that $T_1$ is either essential or there is an essential type (2) annulus.

If $T_1$ is boundary-parallel to the side containing $K$, then it must be boundary-parallel to $\partial N(K)$. But then $A$ would have been boundary-parallel to $\partial N(K)$, a contradiction.

If $T_1$ is boundary-parallel to the side not containing $K$, then the manifold outside $T_1$ is $T \times I$ with $T \times \{0\} = T_1$ and $T \times \{1\} = \partial (F \times I).$  Then $F$ must be an annulus. Since $T_1$ can be isotoped to $\partial (F \times I)$, there is an essential annulus that has one boundary that is on $\partial N(K)$ and a second in $\partial (F \times I)$. This is an essential type (2) annulus. 

The torus $T_1$ is  incompressible to the side containing $K$, as $K$ prevents a compressing disk other than one with boundary parallel to that of $\partial A$, and that cannot exist by incompressibility of $A$. If there is a compressing disk to the other side of $T_1$, compressing along it yields a sphere, which must bound a ball to that side, leaving no place for $\partial (F \times I)$. 

Similar arguments work in the case the two boundaries of $A$  are on the  boundaries of regular neighborhoods of two different components of $L$. 
So we can restrict to the case of an essential  type (2) annulus. 

Suppose that $A$  is such an annulus  with one boundary on $\partial (F \times I)$ and one on $\partial N(K)$ for some component $K$ of $L$. Let $A'$ be the boundary of $N(A \cup K)$. It is now a type (1) annulus. That it is incompressible is immediate from the fact $A$ is incompressible. And similar to the previous cases, it cannot be boundary-parallel to the side containing $K$, and it cannot be boundary-parallel to the other side, because of the presence of $\partial (F \times I)$ to that side. 

Thus, whenever there is an essential type (2) or (3) annulus present, there is an essential torus or an essential type (1)  annulus present.
\end{proof}

We have already eliminated essential tori in Proposition \ref{prop-no-essential-spheres-tori}. So, we now eliminate essential type (1)  annuli. But, for this case, we need all four of the conditions from Theorem \ref{thm-thickened}. 

\begin{lemma} \label{lemma-no-type1-annuli}
Let $M = (F \times I) \setminus N(L)$ as in Theorem \ref{thm-thickened}. Then $M$ has no essential type (1)  annuli. 
\end{lemma}

\begin{proof}
Suppose $\Sigma$ is an essential type (1)  annulus in $M$, and suppose it is meridionally compressible. Then a meridional compression of $\Sigma$ yields two disks with boundary on $\partial (F \times I)$, each punctured once by the link. Let $\Sigma'$ be one of these once-punctured disks; similar to the proof of Proposition \ref{prop-no-essential-disks}, consider the projection of the intersection curves in $\Sigma' \cap F_\pm$ and saddle disks to $\Sigma'$, which is homeomorphic to an annulus. We call the boundary component in $\partial(F \times I)$ the \textit{outer boundary} and we call the boundary component in $\partial N(L)$ the \textit{inner boundary}. After shrinking down the saddle disks, we obtain a 4-valent graph, possibly with some circles without vertices. We denote their union by $\Gamma$ which we call the intersection graph on $\Sigma'$. 

Since $\Sigma$ is punctured exactly once by $L$, it is punctured away from the crossings. In particular, $L$ lies in $F$ near where it punctures $\Sigma$. Furthermore, $\Sigma'$ was obtained via a meridional compression on $\Sigma$, so the core curve of $\Sigma'$ is isotopic in $M$ to a meridian $\mu$ of $L$. Hence, there are exactly two intersection arcs of $\Gamma$ with an endpoint on the inner boundary; all other intersection arcs must have both endpoints on the outer boundary, and in particular, they are trivial on $\Sigma'$.

Observe that we may assume the intersection graph $\Gamma$ has no circles without vertices. Suppose $\gamma$ is such a circle in $\Sigma' \cap F_\pm$. Then it is isotopic via $\Sigma'$ to $\mu$. We know $\gamma$ cannot be trivial in $F_\pm$; otherwise $\mu$ would be trivial in $M$. But then $\gamma$ represents a nontrivial element in $\pi_1(F \times I) \subset \pi_1(M)$ and hence, cannot be homotopic to $\mu$. In particular, this shows that $\Gamma$ is a genuine 4-valent graph on $\Sigma'$, and every complementary region of $\Sigma' \setminus \Gamma$ intersects $\partial \Sigma'$ in its boundary.

Next, we determine when $\Gamma$ contains a fork. As in the proof of Proposition \ref{prop-no-essential-disks} and following remark, the existence of a fork with both edges ending on the outer boundary immediately gives us the desired contradiction by condition (iii). 

Consider the subgraph $\Gamma'$ obtained by throwing away any edges of $\Gamma$ which meet either boundary component of $\Sigma'$. If there are no cycles in $\Gamma'$, then $\Gamma'$ is a collection of trees. 
If one of those trees has two or more vertices, then there are at least two vertices of the tree which appear at the end of leaves. Since in $\Gamma$, every vertex is 4-valent,  each of these vertices has at least three edges ending on $\partial \Sigma'$. Since there are exactly two edges in $\Gamma$ with endpoint on the inner boundary, we know there is at least one fork with both endpoints on the outer boundary. 

If there is only one tree with exactly one vertex, then all four of its edges have an endpoint on $\partial \Sigma'$, and in particular, exactly two (adjacent) edges have endpoint on the outer boundary. This is a fork.  Finally, if there are no vertices, then the two arcs with an endpoint each on the inner boundary have their other endpoint on the outer boundary. Denote these two arcs by $\alpha$ and $\beta$. On the projection surface $F$, $\alpha \cup \beta$ appears as an arc which has both endpoints on  $\partial F$ and intersects $\pi(L)$ transversely in exactly one point (corresponding to the puncture in $\Sigma)$. See Figure \ref{fig-mini-fork} for a diagram of $\Gamma$ on $\Sigma'$ and the corresponding picture on $F$. But $\alpha$ and $\beta$ lie in adjacent regions of $F \setminus \pi(L)$, both of which contain boundary in $\partial F$. This contradicts condition (iii) from the statement of Theorem \ref{thm-thickened}. 

If the intersection graph on $\Gamma'$ does contain a cycle, then because there are no trivial cycles on $\Gamma'$, the only possibility is a single cycle parallel to the boundaries of $\Gamma'$. Because there cannot be any forks, the cycle must intersect two saddles in order that there are two edges ending on the inner boundary. The opposite edge of each saddle ends on the outer boundary of $\Sigma'$. 

The complement of this graph is four disks, two of which live above $F_+$ and two of which live below $F_-$. Let $D_1^+$ denote the disk that lives above $F_+$ and that has a part of its boundary on the inner boundary of $\Sigma'$. We can realize $\partial D_1^+$ as a simple closed curve $\mu$ on $F_+$ that is crossed once by the link corresponding to the inner boundary and twice more by the link at the saddles. Because $F$ is incompressible in $F \times I$, $\mu$ must be trivial on $F \times I$. But we cannot have three strands of the link entering a disk region on $F$, a contradiction. 

\begin{figure} [htbp]
    \centering
    \includegraphics[scale=.6]{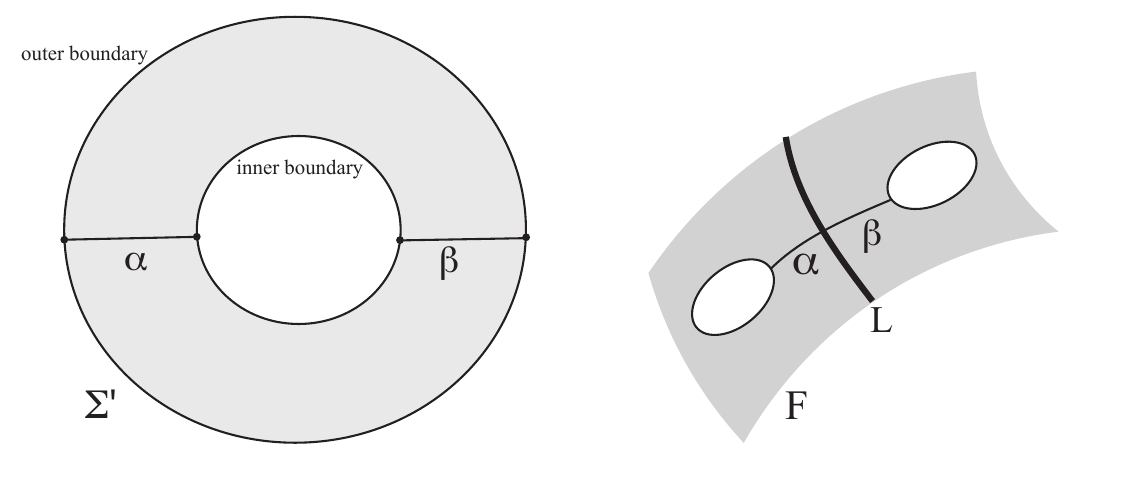}
    \caption{An intersection graph on $\Sigma'$ without vertices and intersection arcs $\alpha$ and $\beta$ as they appear on $F$.}
    \label{fig-mini-fork}
\end{figure}

Thus,  $\Sigma$ must be meridionally incompressible. By Lemma \ref{lemma-at-least-one-int-curve}, without loss of generality $\Sigma \cap F_+$ contains at least one intersection curve $\alpha'$. First, suppose that $\alpha'$ is trivial. If $\alpha'$ is closed, then by Lemma \ref{lemma-triv-int-curve-iff}, it is trivial on $F_+$ and contradicts Lemma \ref{lemma-removing-triv-curve}. Then $\alpha'$ must be an arc (which bounds a disk $D$ in $\Sigma$). Consider the intersection graph $\Gamma'$ on $\Sigma$. By Lemma \ref{lemma-int-curve-one-bubble}, this graph has at least one vertex. 

If $\Gamma'$ has a fork (where the endpoints of its edges can lie on either boundary component of $\Sigma$), then we are done. Suppose that the intersection graph contains a cycle. As the cycle cannot be trivial, it must wrap once around the core curve of the annulus and it must have vertices. To avoid trivial cycles and to avoid forks, the only possibility  is that the remaining two edges coming out of a vertex  must go directly out to the boundary, one to each of the separate boundaries of $\Sigma$.  Thus, we obtain an intersection graph as in Figure \ref{circleinannulus}, but with any even number of vertices. Note that the intersection graph decomposes $\Sigma$ into disks, and any two disks that share an edge in the graph must appear on opposite sides of $F$. This forces the number of vertices to be even.

\begin{figure} [htbp]
    \centering
    \includegraphics[scale=.6]{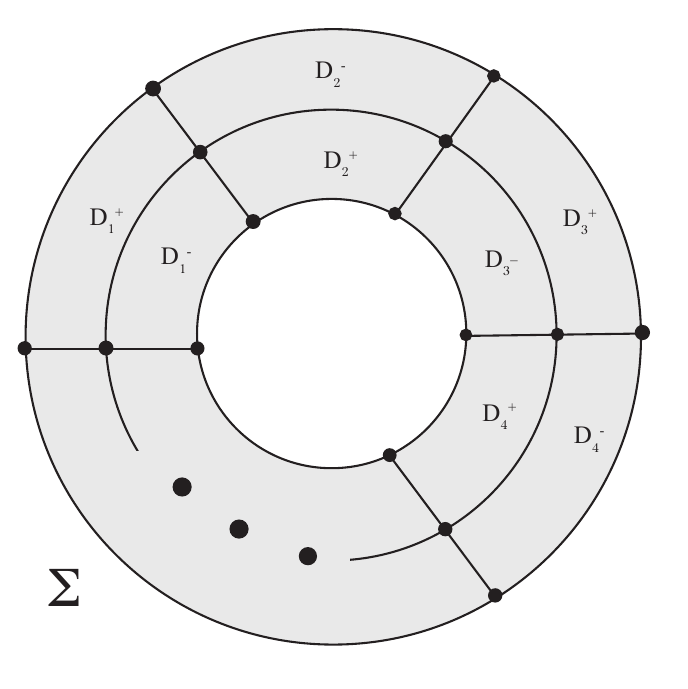}
    \caption{An example of a nontrivial cycle in the intersection graph on $\Sigma$. There must be an even number of vertices on the cycle.}
    \label{circleinannulus}
\end{figure}

We now consider such annuli and show that they contradict condition (iv) of Theorem \ref{thm-thickened}. Suppose there is an annulus with such an intersection graph. Then we first consider how that annulus $\Sigma$ can sit in $F \times I$. Note that if all of the saddles that occur on $\Sigma$ appear in distinct bubbles, then the core curve in the intersection graph yields a simple closed curve $\beta$ on $F$ that passes through the corresponding bubbles, bisecting each crossing. The fact the remaining two edges coming out of each saddle must go directly to a boundary component of $F$ implies that the two complementary regions on $F$ at such a crossing that do not intersect $\beta$ in their interiors must contain these edges and therefore must be annular regions so that there are boundary components of $F$ for these edges to end on. This is exactly the situation that condition (iv) eliminates. Note that in this case, the number of crossings intersected by $\alpha$ is the number of saddles in $\Sigma$, which is even.

Suppose now that there are some saddles from $\Sigma$ that occur in the same bubble. Then the resultant curve $\beta$ on $F$ will pass through a given crossing more than once. However, at any such crossing, $\beta$ must only pass through the same pair of opposite complementary regions, as otherwise all four complementary regions meeting at the crossing would have to be annuli to accommodate the branches on the intersection graph and this contradicts condition (iii). 
The curve $\beta$ must cross itself transversely as it passes through a crossing because it is passing along the diagonal of saddles each time it passes through as in Figure \ref{selfintersectingcurve}. Note that it is still the case that the two opposite complementary regions through which $\beta$ does not pass must still be annuli to accommodate the branches coming out of the saddles. 

\begin{figure} [htbp]
    \centering
    \includegraphics[scale=.4]{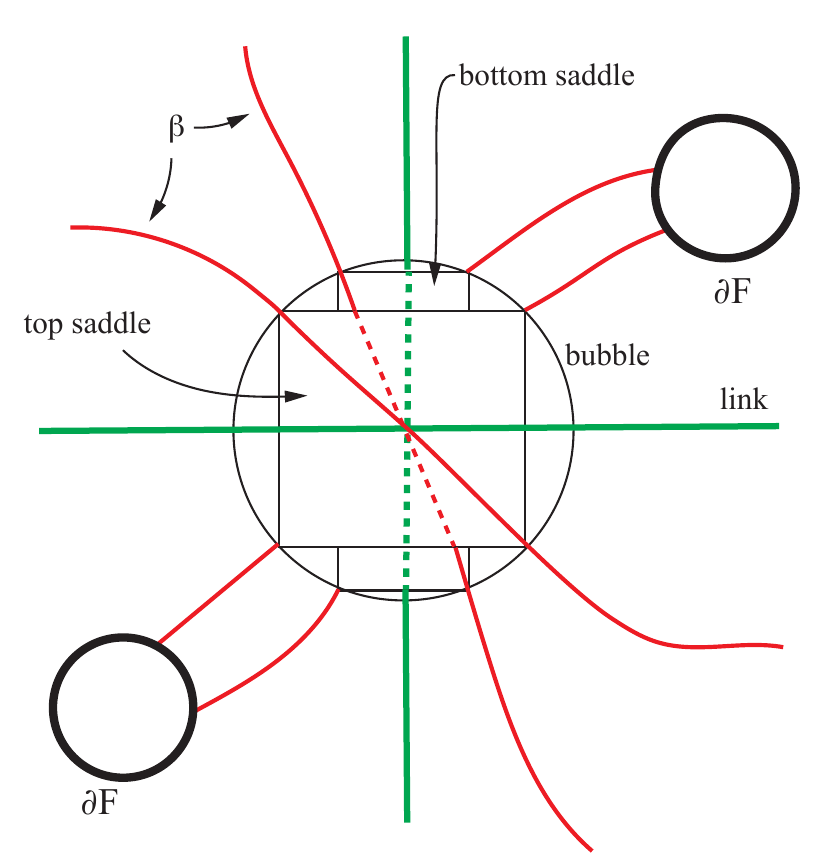}
    \caption{When the curve $\beta$ passes through a crossing more than once, it must intersect the same pair of opposite complementary regions and it must  cross itself transversely in the process.}
    \label{selfintersectingcurve}
\end{figure}

Suppose there are bubbles with more than two saddles. Then we can surger $\beta$ so that it passes through the crossing once if the original passed through an odd number of times and twice if the original passed through an even number of times, as in Figure \ref{surgercurve}. Call the resulting set of curves $\Phi$. If there is at least one component $\phi$ of $\Phi$ with at least one crossing that it passes through once, then surger all remaining crossings that are passed through twice by $\phi$ so as  not to  pass through those crossings. Then take a simple closed curve component $\phi'$ that results and  that passes through at least one crossing. Such a $\phi'$ must exist. It contradicts condition (iv). 

If there is no component of $\Phi$ that passes through a  bubble once, so each component passes through every bubble it intersects twice, then choose any component $\mu$ that passes through bubbles. Its projection on $F$ is the projection of a knot. Because two of the opposite regions at a crossing of $\mu$ must touch regions of the link projection 
that are annuli and the other two regions cannot touch regions of the link projection that are annuli by condition (iii), the projection of $\mu$ must be checkerboard colorable, with the annular regions of the link projection touching a crossing all occurring in the shaded regions of the projection of $\mu$. However, we can then surger the crossings of $\mu$ around the crossings on the boundary of a single shaded region to obtain a loop that passes through each crossing of the link projection once, with annular regions to either side at each crossing. This loop contradicts condition (iv).


\begin{figure} [htbp]
    \centering
    \includegraphics[scale=.4]{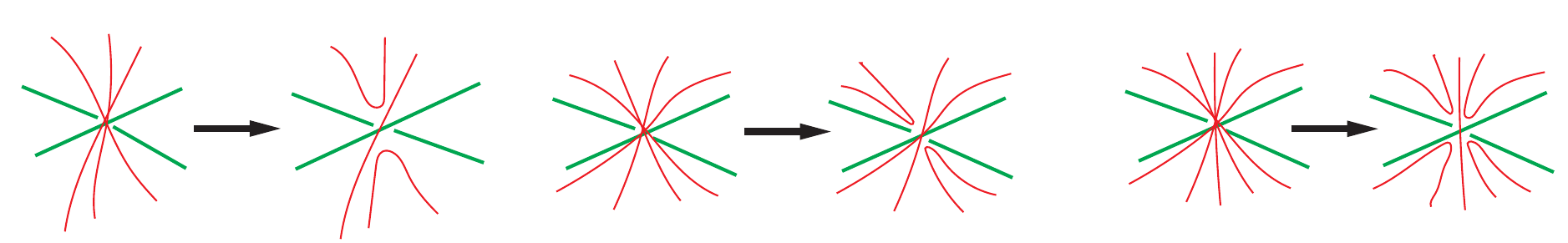}
    \caption{Surgering curve at crossings to make sure it passes through a crossing either once or twice.}
    \label{surgercurve}
\end{figure}

\end{proof}

\begin{proposition} \label{prop-no-essential-annuli}
The manifold $M$ contains no essential annuli. 
\end{proposition}

\begin{proof}
This follows from Proposition \ref{prop-no-essential-spheres-tori} and Lemmas \ref{lemma-type2and3annuliimplytype1} and \ref{lemma-no-type1-annuli}. 
\end{proof}

In order to prove Theorem \ref{thm-thickened} in both directions, we need the following lemma.

\begin{lemma} \label{lemma-failiv-exists-annulus} If there exists a simple closed curve $\alpha$ in $F$ that intersects $\pi(L)$ exactly in a nonempty collection of crossings, such that for each crossing, $\alpha$ bisects the crossing and the two opposite complementary regions meeting at that crossing that do not intersect $\alpha$ near that crossing are annuli, then there exists an essential annulus in $M = (F \times I) \setminus N(L)$. 
\end{lemma}

\begin{figure} [htbp]
    \centering
    \includegraphics[scale=.25]{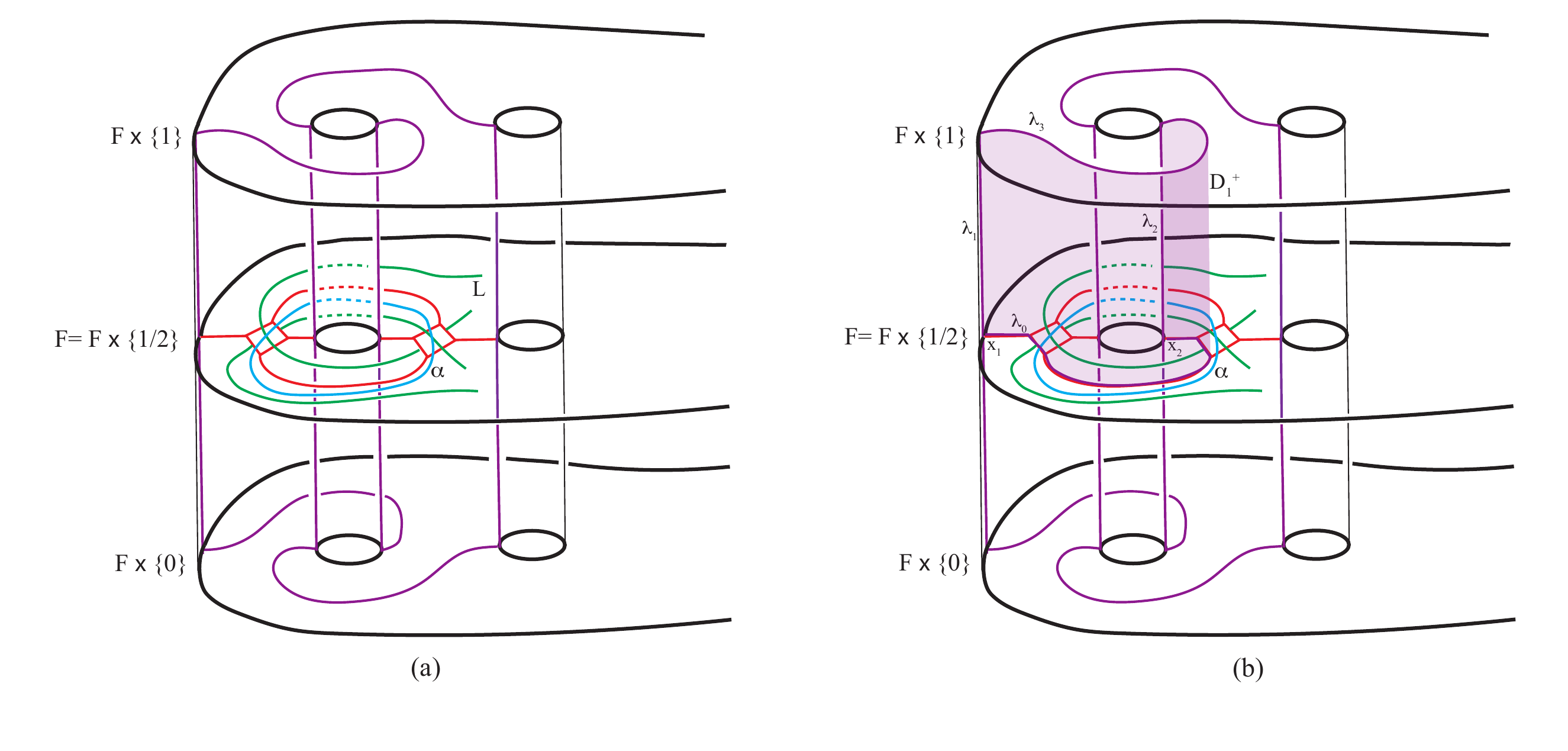}
    \caption{The intersection curves for an annulus that is generated when 
    condition (iv) does not hold.}
    \label{annulusthatfailsiv}
\end{figure}

In Figure \ref{annulusthatfailsiv} (a), we see a particular example of such an annulus. Green curves represent the link, red curves represent the intersection arcs and saddles, and purple curves represent the boundary curves of the annulus.

\begin{proof} First assume that $\alpha$ passes through an even number of crossings $n$. We construct an annulus $\Sigma$ as in Figure \ref{circleinannulus} that exists in the manifold. The nontrivial cycle in $\Sigma$ becomes $\alpha$. Each of the branch edges going out to $\partial F$ on $\Sigma$ exist on $F$ since at each crossing that $\alpha$ passes through, the opposite regions that it does not pass through are annuli. Thus, we can create the corresponding intersection graph on $F$. We describe how to insert the disk $D_1^+$ in $(F \times I) \setminus N(L)$, but the same description works to insert any of the disks $D_i^{\pm}$. The disk $D_1^+$ has boundary on $F = F \times \{1/2\}$ that runs along an arc in the intersection graph starting at $x_1$ on $\partial F$, followed by an arc on the boundary of a saddle, followed by another arc in the intersection graph along $\alpha$, followed by an arc on the boundary of another saddle, followed by an arc in the intersection graph that ends at $x_2$ on $\partial F$ on the opposite side of $\alpha$ from $x_1$. (See Figure \ref{annulusthatfailsiv}(b).) Call this longer arc $\lambda_0.$ The remaining portion of $\partial D_1^+$ lies on the boundary of the handlebody. We describe it as follows. Let $\lambda_1 = \{x_1\} \times [1/2,1]$ and $\lambda_2 = \{x_2\} \times [1/2,1]$. Let $\lambda_3$ be a copy of the arc $\lambda _0$ from $F$ but appearing on $F \times \{1\}$. Then $\partial D_1^+ = \lambda_0 \cup \lambda_1 \cup \lambda_2 \cup \lambda_3$. That this curve bounds a disk that avoids $L$ is immediate from the construction. That the set of disks inserted in this manner together with the corresponding saddles form a properly embedded annulus is also immediate from the construction and the way the various disks share certain edges. It remains to prove that $\Sigma$ is essential.

Suppose $\Sigma$ is compressible. Then $\alpha$ must bound a disk in $F \times I$. But by assumption, since branch arcs to either side of $\alpha$ must end on boundary components of $F$, $\alpha$ cannot be trivial on $F$. Hence, it cannot be trivial in $F \times I$.

Suppose now that $\Sigma$ is boundary compressible. Then if we take a path through the intersection graph that corresponds to two opposite edges at a saddle that each go out to the boundary, together with a diagonal of the saddle, and call that arc $\kappa$, then $\kappa$ together with an arc $\kappa'$ on the boundary of $F \times I$ bound a disk $D$ in $M$. 

Consider the intersection graph of $D$ with $F$.  If $D$ only intersects $F$ along $\kappa$, then $\partial F$ wraps once around $L$ when $\kappa$ passes under $L$ on the saddle, Forcing $D$ to be punctured. Otherwise, there are arcs of intersection that go out to the boundary of $D$. The result is an intersection graph on $D$ much like the intersection graph we had on $D$ when we were proving there were no essential disks in $M$. By the same argument given there using forks, we prove no such disk can exist.

Finally we must consider the case that $\alpha$ passes through an odd number of bubbles. In this case, the same construction yields a properly embedded M\"obius band. The boundary of a regular neighborhood of this M\"obius band is a properly embedded annulus $\Sigma$. It is both incompressible and boundary incompressible to the side that contains the M\"obius band. A similar argument to the previous case demonstrates that it is also incompressible and boundary incompressible to the other side.

\end{proof}

This allows us to complete the proof of Theorem \ref{thm-thickened}, assuming Proposition \ref{prop-weakly-prime}. 

\begin{proof} [Proof of Theorem \ref{thm-thickened}]
By Thurston's Hyperbolization Theorem, $M = \\ (F \times I) \setminus N(L)$ is tg-hyperbolic if and only if there are no essential spheres, tori, disks, or annuli. Then the theorem follows in one direction from Propositions \ref{prop-no-essential-spheres-tori} and \ref{prop-no-essential-disks} and Lemmas \ref{lemma-type2and3annuliimplytype1} and \ref{lemma-no-type1-annuli}, and in the other direction by the comments subsequent to the theorem statement and Lemma \ref{lemma-failiv-exists-annulus}.
\end{proof}

We conclude the section with a proof of Proposition \ref{prop-weakly-prime}. 

\begin{proof} [Proof of Proposition \ref{prop-weakly-prime}]
First, suppose that $\pi(L)$ is not weakly prime, and let $\gamma$ be a circle intersecting $\pi(L)$ transversely in exactly two points such that it bounds a disk $D$ and $D$ contains at least two crossings (note this uses that $\pi(L)$ is reduced). Let $B$ be a regular neighborhood of $D$ which is a 3-ball. Its boundary $\partial B$ is a 2-sphere punctured twice by $L$. Furthermore, $D$ contains at least two (non-reducible) crossings and $\pi(L)$ is alternating; hence, $B$ intersects $L$ in a nontrivial arc and $L$ is not prime. 

It remains to show that if $L$ is not prime, then $\pi(L)$ is not weakly prime. Let $\Sigma \subset F \times I$ be an essential sphere which is punctured twice by $L$. Note that we may assume $\Sigma$ is meridionally incompressible. Indeed, a meridional compression yields two essential spheres, each of which are punctured twice by $L$. Then iteratively perform these compressions until we obtain a meridionally incompressible essential twice-punctured sphere. 

Let $\Tilde{F}$ be the closed orientable surface of genus $g$ obtained by capping off each circle boundary of $F$ with a disk. Then we may regard $L \subset F \times I$ as a link in $\Tilde{F} \times I$ with projection diagram $\pi(L)$ onto $\Tilde{F} \times \{1/2\}$. We refer to this projection surface by $\Tilde{F}$ when there is no ambiguity. Analogously, we may define $\Tilde{F}_+$ (resp. $\Tilde{F}_-$) to be the surfaces obtained from $\Tilde{F}$ by removing the disks in its intersection with the bubbles and replacing them with the upper (resp. lower) hemispheres of the bubbles. Also, we view $\Sigma$ as a sphere in $\Tilde{F} \times I$ which is twice-punctured by $L$ via the inclusion map $F \times I \hookrightarrow \Tilde{F} \times I$. Note that $\Sigma$ is still essential in $\Tilde{F} \times I$.

In Lemma 13 of \cite{small18}, the authors show that when $\Sigma$ is a meridionally incompressible essential sphere which is punctured twice by $L$ and the embedding of $\Sigma$ minimizes $(s, i)$, then there is exactly one intersection curve $\alpha$ in $\Sigma \cap \Tilde{F}_+$ and it intersects $L$ at least twice. This is true even when $\pi(L)$ is not reduced when viewed on $\Tilde{F}_+$. Moreover, the authors of \cite{small18} show that $\alpha$ must be trivial on $\Sigma$ and $\Tilde{F}_+$, and it does not intersect any bubbles. Hence, $\alpha$ is a circle bounding a disk $D$ in $\Tilde{F}_+$ which intersects $\pi(L)$ exactly twice. Since $\Sigma$ is essential, $D$ contains at least one crossing. 

If $\pi(L)$ is reduced when viewed on $\Tilde{F}$, then $\pi(L)$ is not weakly prime on $\Tilde{F}$. Observe that $\alpha$ must bound a disk in $F_+$ as well; otherwise, we find a compression disk for $F \times \{1/2\}$ in $F \times [1/2, 1]$. Hence, $\pi(L)$ is not weakly prime on $F$.

If $\pi(L)$ is not reduced, suppose $\alpha$ may be isotoped slightly in $\Tilde{F}_+$ so it intersects $\pi(L)$ exactly once at a double point. Consider the regions of $\Tilde{F} \setminus \pi(L)$ which are contained in $D$ and do not intersect $\alpha$, and observe that since $\pi(L)$ is reduced on $F$, at least one of these regions is obtained by capping off a boundary component of the corresponding region in $F \setminus \pi(L)$. But then we can find a compression disk for $F \times \{1/2\}$ in $F \times [1/2, 1]$. Hence, the crossings contained in $D$ are unaffected by reducing $\pi(L)$ in $\Tilde{F}_+$, and by the same argument as before, $\pi(L)$ is not weakly prime in $F$. 
\end{proof}

\label{section-pf-thm-thickened} 

\section{Proof of Theorem \ref{surfaceinmanifold}}

\begin{proof}[Proof of Theorem \ref{surfaceinmanifold}] 
Note that $F$ is neither a disk by $\partial$-irreducibility of $Y$ nor an annulus since if it were, we could push a copy off itself and we would have an essential annulus in $Y$ that did not intersect $F$. 

Suppose $\Sigma$ is a properly embedded essential disk or sphere in $Y \setminus N(L)$. If $\Sigma$ does not intersect $F \times I$, its existence contradicts the $\partial$-irreducibility or irreducibility of $Y$. If $\Sigma$ is entirely contained in $(F \times I) \setminus N(L)$ then if it is a sphere, it must bound a ball in $(F \times I) \setminus N(L)$ by Proposition \ref{prop-no-essential-spheres-tori}, which did not use conditions (iii) or (iv), and hence a ball in $Y \setminus N(L)$, a contradiction to its being essential. If it is a disk $D$, the boundary of the disk would have to be a nontrivial curve in  $\partial F \times I$, contradicting the fact $\partial F \times I$ is incompressible in $F \times I$ and also the $\partial$-irreducibility of $Y$.

Still assuming $\Sigma$ is an essential disk or sphere, if $\Sigma$ intersects $F \times I$ but is not entirely contained in it, then by 
incompressibility and $\partial$-incompressibility of $F$, we can replace it by a properly embedded essential $\Sigma'$ that does not intersect $F \times \{0,1\}$,  a contradiction to the cases we have already discussed.

Suppose now that $\Sigma$ is an essential torus or  annulus  in $Y \setminus N(L)$. If it does not intersect $F \times I$, then it must either be compressible or boundary-parallel in $Y$. If it compresses in $Y$, then a compressing disk $D$ must intersect $F \times I$. But by incompressibility of $F$, we can find another compressing disk that does not intersect $F \times I$, contradicting essentiality in $Y \setminus N(L)$. 

If $\Sigma$ is boundary-parallel in $Y$, then there is a boundary component $H$ of $Y$ with which $\Sigma$ is boundary-parallel. In the case $\Sigma$ is a torus, there is a $T \times I $ through which $\Sigma$ is parallel to the boundary, and $H$ is a torus. Since $F$ does not intersect $\Sigma$, it must be contained in $T \times I$ so that $L \subset F \times I$ can prevent $\Sigma$ from being boundary-parallel in $Y \setminus N(L)$. Further, $F$ must have all of its boundary on $H$. But there are no essential surfaces in $T \times I$ with all boundaries on $T \times \{1\}$. So $F$ does not intersect $T \times I$ and $\Sigma$ remains boundary-parallel in $Y \setminus N(L)$, a contradiction.

 In the case $\Sigma$ is an annulus that is boundary-parallel in $Y$, there is a solid torus through which $\Sigma$ is parallel into a boundary component $H$ of $Y$. Then $\Sigma$ must have both boundary components on $H$. Again, there are no essential surfaces to play the role of $F$ in the solid torus with boundary just on $H$, so $F$ does not intersect the solid torus and $\Sigma$ remains boundary-parallel in $Y \setminus N(L)$, a contradiction.

 In the case that $\Sigma$ is entirely contained in $(F \times I) \setminus N(L)$,  $\Sigma$ cannot be a torus by Proposition \ref{prop-no-essential-spheres-tori}, which did not assume conditions (iii) and (iv). 
 
 If $\Sigma$ is an annulus entirely contained in $(F \times I) \setminus N(L)$, then by Lemma \ref{lemma-type2and3annuliimplytype1}, there exists a type (1) annulus $\Sigma'$ that has both boundaries in $\partial F \times I$, 
 and they must be nontrivial curves so that $\Sigma$ is incompressible in $Y \setminus N(L)$. If each boundary is on a different component of $\partial F \times I$, then the two components would be isotopic in $F \times I$  through the annulus, a contradiction to the fact $F$ is not itself an annulus. 
 
 If both boundary components are in the same component of $\partial F \times I$, then
  $\Sigma$ can be isotoped so that $\partial \Sigma$ does not intersect $F$.  Therefore, the intersection graph on $\Sigma$ has no edges that touch $\partial \Sigma$. 
 
 Since the intersection graph is 4-valent, this implies that there are simple closed curves in the intersection graph that are trivial, which contradicts the comment following the proof of Lemma \ref{lemma-triv-int-curve-iff}. Hence the intersection graph is empty. Thus,  $\Sigma$ does not intersect $F$ and it must be boundary-parallel in $F \times I$, a contradiction to its being essential in $Y \setminus N(L)$.



 Let $\Sigma$ be an essential torus or annulus in $Y \setminus N(L)$ that does intersect $F \times I$ but is not entirely contained in it. 
 In the case $\Sigma$ is a torus, we can assume it is meridionally incompressible, as if not, we would generate a twice-punctured sphere that is essential, and Lemma \ref{prop-weakly-prime-circ} would then imply the projection of our link $L$ to $F$ is not weakly prime, contradicting condition (i) of the theorem. 
 
      In the case $\Sigma$ is an annulus, if it is not meridionally incompressible, we can compress to obtain two annuli, each essential. If either had a second meridional compression, we would similarly contradict the fact that the projection of our link $L$ to $F$ is not weakly prime. Thus by replacing $\Sigma$ by one of the two resulting annuli if needed, we can assume $\Sigma$ is both essential and meridionally incompressible. 
      
      We then consider its intersection curves with $F \times \{0, 1\}$. Assume we have chosen $\Sigma$ to minimize the number of such intersection curves. 
     
 Any simple closed intersection curve that bounds a disk on $\Sigma$ must also bound a disk on $F$ by incompressibility. Using irreducibility of $Y$ and of $(F \times I) \setminus N(L)$, we can isotope to eliminate all such simple closed curves, contradicting our choice of $\Sigma$ as having a minimal number of intersection curves. So the only possible simple closed curves of intersection are nontrivial on $\Sigma$ and therefore cut it into annuli. 
     
 In the case when $\Sigma$ is an annulus, we can also have intersection arcs. If such an arc cuts a disk from $\Sigma$, take an outermost such. It must also cut a disk from $F$ by $\partial$-incompressibility of $F$.  The union of these two disks generates a properly embedded disk in $Y$, which by $\partial$-irreducibility of $Y$ must have trivial boundary on $\partial Y$. We can then form a sphere that bounds a ball, through which we can isotope $\Sigma$ to eliminate the intersection arc, again contradicting minimality of intersection curves on $\Sigma$.
     
 In the case there are simple closed curves of intersection, $\Sigma$ will intersect $F \times I$ in a collection of annuli, each  with boundary components on $F \times \{0,1\}$ and possibly one on $\partial N(L)$ and also intersect  $Y \setminus (F \times I)$ in another collection of annuli, the boundaries of which are on $F \times \{0,1\}$. Given such an annulus $A'$ in $(F \times I) \setminus N(L)$ with both boundaries in $F \times \{0, 1\}$, it is incompressible because $\Sigma$ is incompressible in $Y \setminus N(L)$. It cannot be boundary-parallel since we cannot reduce the number of intersections of $\Sigma$ with $F \times \{0, 1\}$. So $A'$ is essential in $(F \times I) \setminus N(L)$. 
 
 But $\partial A'$ does not intersect $F$, and therefore the intersection graph on $A '$ has no edges that touch $\partial A'$. 
 As argued in a previous case, since the intersection graph is 4-valent, this implies that there are simple closed curves in the intersection graph that are trivial, which contradicts the comment following the proof of Lemma \ref{lemma-triv-int-curve-iff}. Hence the intersection graph is empty and  $A'$ does not intersect $F$.  So, it must be boundary-parallel in $F \times I$, a contradiction to the minimality of the number of intersection curves with $F \times \{0, 1\}$.

     

 
The last possibility for $A'$ when the intersection curves of $\Sigma \cap (F \times \{0, 1\})$ are simple closed curves is that  $A'$ intersects $F \times \{0,1\}$ in a single simple closed curve. Then its other boundary is on $\partial N(K)$ where $K$ is a component of $L$. The boundary of a regular neighborhood of $ A' \cup K$ is an annulus $A''$ with both boundaries on one component of $F \times \{0,1\}.$ It must be incompressible since $\Sigma$ is incompressible. 
And it is not boundary-parallel in $(F \times I) \setminus N(L)$ to the side containing $K$. If it was boundary-parallel to the other side, 
 then that side must be a solid torus. The side containing $K$ is already a solid torus. These two solid tori share an annulus on their boundary and their union is all of $F \times I$. Thus, the boundary of $F \times I$ must be a torus and $F$ must be an annulus, a contradiction to our comment at the beginning of the proof that $F$ cannot be an annulus. 

Thus, $A''$ is essential in $F \times I$ with both boundaries on $F \times \{0,1\}$. But we have already eliminated this possibility. 





Finally, suppose that $\Sigma$ is an annulus, and all intersection curves are arcs from one boundary of $\Sigma$ to the other. They cut $\Sigma$ into disks that alternately lie in $(F \times I) \setminus N(L)$ and $Y \setminus (F \times I)$. Let $D$ be such a disk in $(F \times I) \setminus N(L)$. Its boundary $\partial D$ breaks up into four arcs, two of which are non-adjacent arcs in $F \times \{0, 1\}$ and two of which are in $\partial F \times I$.

The disk $D$ must intersect $F \times \{1/2\}$ or we could isotope it off $F \times I$ and reduce the number of intersections of $\Sigma$ with  $F \times I$. Then we can isotope to obtain an intersection graph on $D$ with at least one vertex by Lemma \ref{lemma-int-curve-one-bubble}. Hence as argued previously, there must be a fork. 

In fact, we can argue that there are at least four forks. If there is just a single vertex, there are four edges that go to the boundary and hence four forks. If there is more than one vertex, removing all edges in the intersection graph that intersect $\partial D$, we are left with a tree or trees with two or more leaves. Each then has three edges going out to the boundary, and thus generates at least two forks. So we have at least four forks. 

If both edges of a fork end on the same boundary component of $F$, then we create a region in the complement of the intersection graph on $F$ that has one strand of $L$ entering from a bubble with nowhere to go, since the region is a disk or annulus, a contradiction. Hence, the two edges making up the fork must end on distinct components of the boundary of $F$. However,
in order for $\partial D$ to pass from one component of $\partial F \times I$ to a different component, it must travel on the boundary of $F \times I$ up $\partial F \times I$, across one of $F \times \{0, 1\}$ and down $\partial F \times I$. Since there are at least four such forks, this contradicts the fact there are only two connected arcs on the boundary of $D$ that are in $\partial F \times I$. 


Therefore, since there are no essential spheres, disks, tori or annuli in $Y \setminus N(L)$, it must be tg-hyperbolic.
    
\end{proof}

We conclude this section with a proof of Proposition \ref{prop-weakly-prime-circ}. 

\begin{proof} [Proof of Proposition \ref{prop-weakly-prime-circ}]
First, the direction ``$\pi(L)$ is not weakly prime implies $L$ is not prime in $Y$" follows as in the proof of Proposition \ref{prop-weakly-prime}. 
To prove the converse, suppose that $L$ is not prime and let $\Sigma \subset Y$ be an essential sphere punctured twice by $L$. Let $G = F \times \{0\} \cup F \times \{1\}$. Consider the intersection curves in $\Sigma \cap G$. If the intersection is empty, then we can cut along $G$ so that $\Sigma$ is an essential sphere embedded in $F \times I$. Then by Proposition \ref{prop-weakly-prime}, $\pi(L)$ is not weakly prime on $F$. 

Now suppose $\Sigma \cap G$ is not empty, and let $\alpha \subset \Sigma \cap G$ be an intersection curve which is innermost on $\Sigma$. Since $G$ is incompressible in $Y$, $\alpha$ bounds a disk $D \subset G$ on $G$. Let $D' \subset \Sigma$ be a disk bounded by $\alpha$ on $\Sigma$ which contains no other intersection curves of $\Sigma \cap G$: then $D \cup D'$ is a 2-sphere in $Y \setminus G$.
If $D \cup D'$ is outside $F \times I$ it must bound a ball outside $F \times I$ and we can isotope to remove the intersection. If it is inside $F \times I$, it bounds a ball in $F \times I$ by Proposition \ref{prop-no-essential-spheres-tori}.

Thus we can eliminate all intersections of $\Sigma $ with $G$, and hence Proposition \ref{prop-weakly-prime} implies  $\pi(L)$ is not weakly prime.
\end{proof}

\label{section-pf-thm-circle-bundle}  


\section{Applications and Further Directions}
One motivation for Theorem \ref{thm-thickened} comes from studying hyperbolic links in handlebodies. In addition to being interesting objects in their own right as a natural generalization of classical link theory, they also show up naturally when studying hyperbolicity of knotoids and generalized knotoids as in \cite{generalizedknotoids} and \cite{hypknotoids}. Indeed, the map $\phi_{\Sigma}^D$ constructed in these papers allows us to associate a hyperbolic volume to generalized knotoids by mapping them to the set of links in a handlebody. 

Observe that a genus $g$ handlebody can be obtained by thickening a 2-sphere with $(g + 1)$ disks removed, or more generally by thickening a genus $k$ closed orientable surface with $(g - k)$ disks removed (where $k \geq 1$). In the remainder of this section, when we refer to a projection surface in a handlebody, we always mean one of these surfaces, so that the handlebody is obtained by thickening it. Then Theorem \ref{thm-thickened} is useful for studying links in handlebodies, and hence, for studying generalized knotoids as well. 

As one application of the theorem, note that if a generalized knotoid $k$ has any poles of nonzero valency, then $\phi_\Sigma^D$ never yields a link which is cellular alternating with respect to one of these projection surfaces. This is because the construction requires us to double the rail diagram of $k$ across the boundary portions corresponding to the poles of nonzero valency. 

However, we can restrict to the class of generalized knotoids whose poles are all valency-zero, that is, generalized knotoids whose diagram consists of a link on $\Sigma$ together with a set of valency-zero poles. These are the \textit{staked links} defined in \cite{generalizedknotoids}. Then, as noted in Proposition 7.7 of \cite{generalizedknotoids}, Theorem \ref{thm-thickened} precisely characterizes which alternating staked links (or equivalently, alternating links in handlebodies) are hyperbolic under $\phi_\Sigma^D$.

In that paper, this is used to prove Theorem 7.8, which says that every link with a checkerboard-colorable diagram on a closed surface $\Sigma$ has a diagram such that staking that diagram makes the resulting link hyperbolic in $\Sigma \times I$. In particular, this means that we can define the staked volume for any such link to be the minimum volume of any hyperbolic staking of the link. Since all link diagrams in $S^3$ are checkerboard-colorable, we can define staked volume for every link in $S^3$. See the last section of \cite{generalizedknotoids} for more details.

\begin{figure} [htbp]
    \centering
    \includegraphics[scale = 0.7]{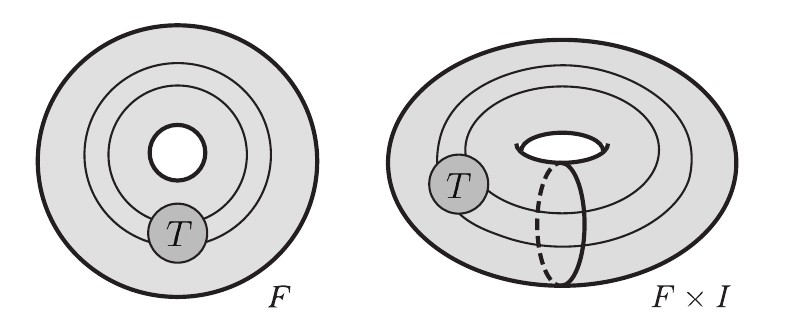}
    \caption{Here $F$ is an annulus and $T$ is a prime, cellular alternating tangle which is not an integer tangle. By Theorem \ref{thm-thickened}, this is a family of tg-hyperbolic links in a solid torus. There are no closed projection surfaces which satisfy the hypotheses of Theorem 1.1 from \cite{hp17}. For example, if we choose the torus parallel to $\partial(F \times I)$, the link will not have a cellular projection.} \label{fig-thm1.1-nonex}
\end{figure}

As we mentioned in the introduction, there are examples of links in handlebodies shown to be hyperbolic by Theorem \ref{thm-thickened} that are not covered by the hypothesis of Theorem 1.1 from \cite{hp17}. For example, consider the family of examples in Figure \ref{fig-thm1.1-nonex}.

We would like to extend Theorem \ref{thm-thickened} to cellular non-alternating links. The following result gives an extension in this direction:

\begin{corollary}
Let $F$ be a projection surface with nonempty boundary which is not a disk, and let $L \subset F \times I$ be a link with a reduced cellular (not necessarily alternating) projection diagram $\pi(L) \subset F \times \{1/2\}$, and let $M = (F \times I) \setminus N(L)$. Suppose conditions (i)-(iv) of Theorem \ref{thm-thickened} are satisfied, as well as the following:
\begin{enumerate} [label = (\roman*)]
\setcounter{enumi}{4}
    \item let $c_1, \dotsc, c_n$ be crossings of $\pi(L)$ such that $\pi(L)$ becomes alternating after each $c_i$ is changed to the opposite crossing. Each $c_i$ locally divides $F$ into four complementary regions such that a pair of opposite regions are homeomorphically annuli. 
\end{enumerate}
Then the conclusion of Theorem \ref{thm-thickened} holds. 
\end{corollary}


\begin{proof}
Consider the projection diagram $\pi(L) \subset F$. For each crossing $c$ in $\{c_i\}_{i = 1}^n$, let $R_1$ and $R_2$ denote the two complementary regions of $F \setminus \pi(L)$ which meet $c$ and are homeomorphically annuli. There is an arc $\alpha \subset F$ which has an endpoint each on $\partial R_1 \cap \partial F$ and $\partial R_2 \cap \partial F$ and intersects $\pi(L)$ exactly once through $c$. Then $\alpha \times I$ is a properly embedded disk $\Sigma$ in $F \times I$ which is punctured twice by $L$. See Figure \ref{fig-crossing-change}. 

\begin{figure}[htbp]
    \centering
    \includegraphics[scale = 0.85]{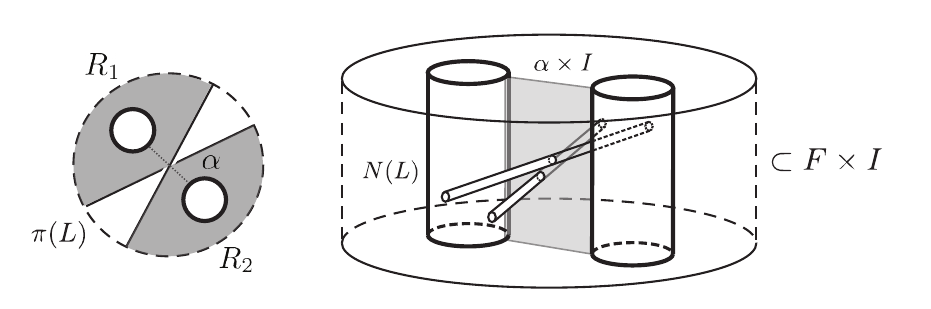}
    \caption{On the left is a local picture of $\pi(L) \subset F$ near crossing $c$ with the $R_i$ shaded. Crossing the arc $\alpha$ by $I$ yields the twice-punctured disk $\Sigma$, shown on the right.}
    \label{fig-crossing-change}
\end{figure}
Now we may cut $M$ along $\Sigma$, yielding two copies, $\Sigma_1$ and $\Sigma_2$. Reglue $\Sigma_2$ to $\Sigma_1$ along a rotation by $2\pi$: this has the effect of changing the crossing $c$ to the opposite crossing, and the resulting manifold is homeomorphic to $M$. By performing this operation for each of the $c_i$, we obtain a link $L' \subset F \times I$ such that $(F \times I) \setminus N(L')$ is homeomorphic to $M$ and its projection $\pi(L')$ is cellular alternating on $F$. Moreover, conditions (i), (ii), (iii) and (iv) still hold. Then the statement follows from Theorem \ref{thm-thickened}. 
\end{proof}

The corollary expands the number of known hyperbolic staked links as defined in \cite{generalizedknotoids}, or equivalently tg-hyperbolic links in handlebodies. 
Theorem \ref{thm-thickened} may also be combined with results from \cite{complinks} to give other ways of obtaining tg-hyperbolic links in handlebodies, namely, by \textit{composition}.

There is much work to be done in expanding the number of links in handlebodies known to be tg-hyperbolic. The methods used in this paper rely heavily on the alternating property: however, it is conceivable that these methods might be adapted for almost alternating links by taking into account the different behavior of the intersection curves at the non-alternating crossing.


Another direction is to shrink the number of hypotheses needed for  Theorem \ref{surfaceinmanifold}. Theorem 1.1 of \cite{hp17} is very powerful in this sense: it applies to links in an arbitrary compact 3-manifold $Y$ (satisfying some mild conditions) with a cellular alternating diagram on a \textit{closed} projection surface that is not necessarily incompressible in $Y$ but rather the diagram satisfies a certain representivity condition. There should be a version of Theorem \ref{surfaceinmanifold} where $F$ need not be incompressible and $\partial$-incompressible, but similarly, the diagram satisfies an appropriate representativity condition. 

We might also try to generalize Theorems \ref{thm-thickened} and \ref{surfaceinmanifold} to allow for nonorientable projection surfaces $F$ or for nonorientable $I$-bundles. Since the analogous results for closed surfaces in \cite{small18} hold in the nonorientable case, we suspect these generalizations hold here as well. 

We are also interested in volume computations for hyperbolic alternating links in thickened surfaces with boundary. In \cite{lackenby}, Lackenby proves a lower bound on hyperbolic volume for alternating links in $S^3$ in terms of the number of twist regions, which can be read off the link diagram. Howie and Purcell generalize this in \cite{hp17} to a lower bound for volumes of links in $Y$. It would be interesting to try adapting their methods to prove a similar lower bound on volume in our case. This might be done by defining a slightly more general version of Howie and Purcell's \textit{angled chunks} which can account for boundary coming from the projection surface. Alternatively, we might try to find proofs of the lower bounds from the viewpoint of bubbles instead.

\label{section-apps}

\bibliographystyle{plain}
\nocite{*}
\bibliography{references}

\begin{thebibliography}{10}

\bibitem{adamsnew}
C.~Adams.
\newblock Quasifuchsian surfaces in hyperbolic knots complements.
\newblock {\em J. Austr. Math. Soc.}, 55:116--131, 1993.

\bibitem{adams}
C.~Adams.
\newblock Toroidally alternating knots and links.
\newblock {\em Topology}, 33(2):353--369, 1994.

\bibitem{adams-forks}
C.~Adams.
\newblock Generalized augmented alternating links and hyperbolic volumes.
\newblock {\em Algebraic {\&} Geometric Topology}, 17(6):3375--3397, 2017.

\bibitem{small18}
C.~Adams, C.~{Albors-Riera}, B.~Haddock, Z.~Li, D.~Nishida, B.~Reinoso, and
  L.~Wang.
\newblock Hyperbolicity of links in thickened surfaces.
\newblock {\em Topology and its Applications}, 256(1):262--278, 2019.

\bibitem{generalizedknotoids}
C.~Adams, A.~Bonat, M.~Chande, J.~Chen, M.~Jiang, Z.~Romrell, D.~Santiago,
  B.~Shapiro, and D.~Woodruff.
\newblock Generalizations of knotoids and spatial graphs.
\newblock {\em ArXiv:2209.01922}, 2022.

\bibitem{hypknotoids}
C.~Adams, A.~Bonat, M.~Chande, J.~Chen, M.~Jiang, Z.~Romrell, D.~Santiago,
  B.~Shapiro, and D.~Woodruff.
\newblock Hyperbolic knotoids.
\newblock {\em ArXiv:2209.04556}, 2022.

\bibitem{complinks}
C.~Adams and D.~Santiago.
\newblock Composition properties of hyperbolic links in handlebodies, 2023.

\bibitem{ast}
I.~Agol, P.~Storm, and W.~Thurston.
\newblock Lower bounds on volumes of hyperbolic {H}aken 3-manifolds.
\newblock {\em Journal of the American Mathematical Society}, 20(4):1053--1077,
  2007.

\bibitem{cfw}
D.~Calegari, M.~Freedman, and K.~Walker.
\newblock Positivity of the universal pairing in 3 dimensions.
\newblock {\em Journal of the American Mathematical Society}, 23(1):107–188,
  2010.

\bibitem{FMP}
R.~Frigerio, B.~Martelli, and C.~Petronio.
\newblock Small hyperbolic 3-manifolds with geodesic boundary.
\newblock {\em Experiment. Math.}, 13(2):171--184, 2004.

\bibitem{Frigerio}
Roberto Frigerio.
\newblock An infinite family of hyperbolic graph complements in ${S}^3$.
\newblock {\em J. Knot Theory Ramifications}, 14(4):479–496, 2005.

\bibitem{hatcher}
A.~Hatcher.
\newblock Notes on basic 3-manifold topology, 2007.

\bibitem{hp17}
J.A. Howie and J.~Purcell.
\newblock Geometry of alternating links on surfaces.
\newblock {\em Trans. Amer. Math. Soc.}, 373:2349--2397, 2020.

\bibitem{lackenby}
M.~Lackenby.
\newblock The volume of hyperbolic alternating link complements, 2000.

\bibitem{martelli}
B.~Martelli.
\newblock {\em An Introduction to Geometric Topology}.
\newblock arXiv, 2016.

\bibitem{menasco}
W.~Menasco.
\newblock Closed incompressible surfaces in alternating knot and link
  complements.
\newblock {\em Topology}, 23(1):37--44, 1984.

\bibitem{simplesmallknots}
R.~Qiu and S.~Wang.
\newblock Simple, small knots in handlebodies.
\newblock {\em Topology and its Applications}, 144:211--227, 2004.

\end{thebibliography}
\end{document}